\documentclass[12pt,twoside,british]{article}
\usepackage[T1]{fontenc}
\usepackage[utf8]{inputenc}
\usepackage[a4paper]{geometry}
\usepackage{color}
\usepackage{amsmath}
\usepackage{amsthm}
\usepackage{amssymb}
\usepackage{setspace}

\makeatletter
\numberwithin{equation}{section}
\numberwithin{figure}{section}
\theoremstyle{plain}
\newtheorem{thm}{\protect\theoremname}[section]
\theoremstyle{plain}
\newtheorem{lem}[thm]{\protect\lemmaname}
\theoremstyle{definition}
\newtheorem{defn}[thm]{\protect\definitionname}
\theoremstyle{plain}
\newtheorem{prop}[thm]{\protect\propositionname}
\newcommand{\lyxaddress}[1]{
	\par {\raggedright #1
	\vspace{1.4em}
	\noindent\par}
}

\@ifundefined{date}{}{\date{}}
\usepackage[english]{babel}
\usepackage{latexsym}
\usepackage{amsthm}
\usepackage{eucal}
\usepackage{eufrak}
\usepackage{epstopdf}
\usepackage{graphicx}
\theoremstyle{plain}

\newtheoremstyle{boldremark}
    {\dimexpr\topsep/2\relax} 
    {\dimexpr\topsep/2\relax} 
    {}          
    {}          
    {\bfseries} 
    {.}         
    {.5em}      
    {}          

\theoremstyle{boldremark}
\newtheorem{brem} [thm] {Remark} 

\usepackage{fancyhdr}
\usepackage{blindtext}
\usepackage{titlesec}

\titleformat{name=\chapter}[display]
{\normalfont\Huge\itshape}
{\titlerule[1pt]\vspace{-33pt}\filleft%
  \parbox[t]{6em}{%
    \raggedleft%
    \rule{\linewidth}{0.5ex}\newline%
    \chaptertitlename\ \thechapter%
  }%
}
{4pc}
{\normalfont\upshape\bfseries\Huge}
\makeatother

\usepackage{babel}
\providecommand{\definitionname}{Definition}
\providecommand{\lemmaname}{Lemma}
\providecommand{\propositionname}{Proposition}
\providecommand{\theoremname}{Theorem}

\begin{document}
\title{\noindent \textbf{Regularity results for a class of widely\\degenerate
parabolic equations}}
\author{\noindent Pasquale Ambrosio, Antonia Passarelli di Napoli \thanks{P. Ambrosio and A. Passarelli di Napoli have been partially supported
by the Gruppo Nazionale per l'Analisi Matematica, la Probabilità e
le loro Applicazioni (GNAMPA) of the Istituto Nazionale di Alta Matematica
(INdAM). A. Passarelli di Napoli has been partially supported by Università
degli Studi di Napoli ``Federico II'' through the Project FRA (000022-75-2021-FRA-PASSARELLI).}}
\date{\noindent September 1, 2023}
\maketitle
\begin{abstract}
\begin{singlespace}
\noindent Motivated by applications to gas filtration problems, we
study the regularity of weak solutions to the strongly degenerate
parabolic PDE
\[
u_{t}-\mathrm{div}\left((\vert Du\vert-\nu)_{+}^{p-1}\frac{Du}{\vert Du\vert}\right)=f\,\,\,\,\,\,\,\,\,\mathrm{in}\,\,\,\Omega_{T}=\Omega\times(0,T),
\]
where $\Omega$ is a bounded domain in $\mathbb{R}^{n}$ for $n\geq2$,
$p\geq2$, $\nu$ is a positive constant and $\left(\,\cdot\,\right)_{+}$
stands for the positive part. Assuming that the datum $f$ belongs
to a suitable Lebesgue-Sobolev parabolic space, we establish the Sobolev
spatial regularity of a nonlinear function of the spatial gradient
of the weak solutions, which in turn implies the existence of the
weak time derivative $u_{t}$. The main novelty here is that the structure
function of the above equation satisfies standard growth and ellipticity
conditions only outside a ball with radius $\nu$ centered at the
origin. We would like to point out that the first result obtained
here can be considered, on the one hand, as the parabolic counterpart
of an elliptic result established in \cite{Br}, and on the other
hand as the extension to a strongly degenerate context of some known
results for less degenerate parabolic equations.\vspace{0.2cm}
\end{singlespace}
\end{abstract}
\begin{singlespace}
\end{singlespace}
\begin{singlespace}

\section{Introduction and statement of the results}
\end{singlespace}

\begin{singlespace}
\noindent $\hspace*{1em}$In this paper, we study the local regularity
properties of weak solutions $u:\Omega_{T}\rightarrow\mathbb{R}$
to strongly degenerate parabolic equations of the type
\begin{equation}
u_{t}-\mathrm{div}\left((\vert Du\vert-\nu)_{+}^{p-1}\frac{Du}{\vert Du\vert}\right)=f\,\,\,\,\,\,\mathrm{in}\,\,\,\Omega_{T}=\Omega\times(0,T),\label{eq:1}
\end{equation}

\noindent for exponents $p\geq2$, where $\Omega$ is a bounded domain
in $\mathbb{R}^{n}$ ($n\geq2$), $T>0$, $\nu$ is a positive constant
and $\left(\,\cdot\,\right)_{+}$ stands for the positive part. \\
$\hspace*{1em}$The main feature of this PDE is that the structure
function satisfies standard growth and ellipticity conditions for
a growth rate $p\geq2$, but only outside a ball with radius $\nu$
centered at the origin.\\
$\hspace*{1em}$The elliptic version of the above equation naturally
arises in optimal transport problems with congestion effects, and
the regularity properties of its weak solutions have been widely investigated:
see, for example, \cite{Ambr,Bog,Br0} and \cite{Br}. As far as we
know, no parabolic counterpart of such works is available in the literature.
On the other hand, we would like to point out that a motivation for
studying equations of the type (\ref{eq:1}) can be found in Section
\ref{subsec:Motivation} below.\\
$\hspace*{1em}$Here, we establish the Sobolev spatial regularity
of a nonlinear function of the spatial gradient $Du$ of the weak
solutions to equation (\ref{eq:1}) (see Theorem \ref{thm:main4}
below), which in turn implies the Sobolev time regularity of the solutions
(cf. Theorem \ref{thm:timeregularity}), by assuming that the datum
$f$ belongs to a suitable Lebesgue-Sobolev parabolic space. These
results are obtained by adapting the techniques for the evolutionary
$p\,$-Laplacian to this more degenerate context. In fact, for less
degenerate parabolic problems, these issues have been widely investigated,
as one can see, for example, in \cite{Duzaar,Giann} (where $f=0$)
and in \cite{Sche}. Moreover, establishing the Sobolev regularity
of the solutions with respect to time, once the higher differentiability
in space has been obtained, is a quite usual fact in these problems:
see, for instance, \cite{Lind1,Lind2,Lind3}.\\
$\hspace*{1em}$The distinguishing feature of equation (\ref{eq:1})
is that the principal part behaves like a \textit{$p$-Laplace operator
only at infinity}. Before giving the main results of this paper, let
us summarize a few previous results on this topic: the regularity
of solutions to parabolic problems with asymptotic structure of $p\,$-Laplacian
type has been explored in \cite{Isernia}, where a BMO regularity
has been proved for solutions to asymptotically parabolic systems
in the case $p=2$ and $f=0$ (see also \cite{Kuusi}, where the local
Lipschitz continuity of weak solutions with respect to the spatial
variable is established). In addition, we want to mention the results
contained in \cite{Byun}, where nonhomogeneous parabolic problems
involving a discontinuous nonlinearity and an asymptotic regularity
in divergence form of $p\,$-Laplacian type are considered. There,
the authors establish a global Calderón-Zygmund estimate by converting
a given asymptotically regular problem to a suitable regular problem.\\
$\hspace*{1em}$One of the main novelties of this work is an observation
of interpolative nature, which allows us to suitably weaken the assumptions
on the datum $f$: this comes from an idea that has already been exploited
in the recent paper \cite{Clop}, in the elliptic setting. In fact,
our assumption on the regularity of $f$ is weaker than those considered
in the mentioned works.\\
$\hspace*{1em}$The first result we prove in this paper is the following
theorem, which can be considered as the parabolic counterpart of Theorem
4.2 in \cite{Br}. We refer to Section \ref{sec:prelim} for notation
and definitions.
\end{singlespace}
\begin{thm}
\begin{singlespace}
\noindent \label{thm:main4} Let $n\geq2$, $p\geq2$, $\frac{np\,+\,4}{np\,+\,4\,-n}\leq\vartheta<\infty$
and $f\in L^{\vartheta}\left(0,T;W^{1,\vartheta}(\Omega)\right)$.
Moreover, assume that 
\[
u\in C^{0}\left((0,T);L^{2}(\Omega)\right)\cap L^{p}\left(0,T;W^{1,p}(\Omega)\right)
\]
is a weak solution of equation $\mathrm{(\ref{eq:1})}$. Then the
solution satisfies 
\[
H_{\frac{p}{2}}(Du)\,\in\,L_{loc}^{2}\left(0,T;W_{loc}^{1,2}(\Omega,\mathbb{R}^{n})\right),
\]
where 
\[
H_{\frac{p}{2}}(Du):=\left(\left|Du\right|-\nu\right)_{+}^{p/2}\frac{Du}{\left|Du\right|}.
\]
Furthermore, the following estimate\begin{equation}\label{eq:faible}
\begin{split}
\int_{Q_{\varrho/2}(z_{0})}\left|DH_{\frac{p}{2}}(Du)\right|^{2}dz\,&\leq\,\,c\left(\nu\,\Vert Df\Vert_{L^{\vartheta}(Q_{R_{0}})}+\,\Vert Df\Vert_{L^{\vartheta}(Q_{R_{0}})}^{\frac{np\,+\,4}{np\,+\,2\,-\,n}}\right)\\
&\,\,\,\,\,\,\,+\frac{c}{R^{2}}\left(\Vert Du\Vert_{L^{p}(Q_{R_{0}})}^{p}+\,\Vert Du\Vert_{L^{p}(Q_{R_{0}})}^{2}+\nu^{p}+\nu^{2}\right)
\end{split}
\end{equation}holds true for any parabolic cylinder $Q_{\varrho}(z_{0})\subset Q_{R}(z_{0})\subset Q_{R_{0}}(z_{0})\Subset\Omega_{T}$
and a positive constant $c$ depending at most on $n$, $p$, $\vartheta$
and $R_{0}$.
\end{singlespace}
\end{thm}

\begin{singlespace}
\noindent $\hspace*{1em}$One of the main tools in the proof of Theorem
\ref{thm:main4} is the difference quotients technique used in the
spatial directions. Here we will argue as in \cite[Lemma 5.1]{Duzaar}
and \cite[Theorem 4.1]{Giann}, but we need to take into account the
strong degeneracy of equation (\ref{eq:1}). This is why we obtain
the Sobolev spatial regularity not for the usual function $V_{p}(Du):=\left|Du\right|^{\frac{p-2}{2}}Du$,
but for the vector field $H_{\frac{p}{2}}(Du)$, which vanishes in
the set where equation (\ref{eq:1}) becomes degenerate.\\
The proof of Theorem \ref{thm:main4} is also based on a comparison
argument with the solutions of a family of much less degenerate parabolic
problems having a smooth inhomogeneity. To the solutions of these
problems we can apply the\textit{ a priori }estimates that we establish
in Section \ref{sec:A priori}, whose constants are independent of
the less degenerate principal parts. Thereafter, we show that the
$L^{p}$-norms of the spatial gradients of such solutions are uniformly
bounded, and this allows us to transfer the higher differentiability
in space of the comparison maps to the solution of our equation (see
Section \ref{sec:main theo} below).

\noindent $\hspace*{1em}$As we anticipated earlier, from the previous
result we can easily deduce that $u$ admits a weak time derivative
$u_{t}$, which belongs to the local Lebesgue space $L_{loc}^{\min\,\{\vartheta,\,p'\}}(\Omega_{T})$,
where $p'=p/(p-1)$ is the conjugate exponent of $p$. The idea is
roughly as follows. Consider equation (\ref{eq:1}); since the above
theorem tells us that in a certain pointwise sense the second spatial
derivatives of $u$ exist, then we may develop the expression under
the divergence symbol; this will give us an expression that equals
$u_{t}$, from which we get the desired summability of the time derivative.
Such an argument must be made more rigorous. Furthermore, we also
need to make explicit \textit{a priori} local estimates. These are
provided in the following
\end{singlespace}
\begin{thm}
\begin{singlespace}
\noindent \label{thm:timeregularity} Under the assumptions of Theorem
\ref{thm:main4}, the time derivative of the solution exists in the
weak sense and satisfies
\[
\partial_{t}u\,\in\,L_{loc}^{\min\,\{\vartheta,\,p'\}}(\Omega_{T}).
\]
Furthermore, the following estimate\begin{equation}\label{eq:estcor}
\begin{split}
&\left(\int_{Q_{\varrho/2}(z_{0})}\left|\partial_{t}u\right|^{\min\,\{\vartheta,\,p'\}}\,dz\right)^{\frac{1}{\min\,\{\vartheta,\,p'\}}}\\
&\,\,\,\,\,\,\,\leq\,\,c\,\Vert f\Vert_{L^{\vartheta}(Q_{R_{0}})}+\,c\,\,\Vert Du\Vert_{L^{p}(Q_{R_{0}})}^{\frac{p-2}{2}}\,\left(\nu\,\Vert Df\Vert_{L^{\vartheta}(Q_{R_{0}})}+\,\Vert Df\Vert_{L^{\vartheta}(Q_{R_{0}})}^{\frac{np\,+\,4}{np\,+\,2\,-\,n}}\right)^{\frac{1}{2}}\\
&\,\,\,\,\,\,\,\,\,\,\,\,\,\,+\frac{c}{R}\left(\Vert Du\Vert_{L^{p}(Q_{R_{0}})}^{2p-2}+\,\Vert Du\Vert_{L^{p}(Q_{R_{0}})}^{p}+(\nu^{p}+\nu^{2})\,\Vert Du\Vert_{L^{p}(Q_{R_{0}})}^{p-2}\right)^{\frac{1}{2}}
\end{split}
\end{equation}

\noindent holds true for any parabolic cylinder $Q_{\varrho}(z_{0})\subset Q_{R}(z_{0})\subset Q_{R_{0}}(z_{0})\Subset\Omega_{T}$
and a positive constant $c$ depending on $n$, $p$, $\vartheta$
and $R_{0}$. 
\end{singlespace}
\end{thm}

\begin{singlespace}
\noindent $\hspace*{1em}$To conclude this introduction, it is worth
pointing out that, starting from the weaker assumption $f\in L^{\frac{np\,+\,4}{np\,+\,4\,-n}}\left(0,T;W^{1,\frac{np\,+\,4}{np\,+\,4\,-n}}(\Omega)\right)$,
Sobolev regularity results such as those of Theorems \ref{thm:main4}
and \ref{thm:timeregularity} seem not to have been established yet
for weak solutions to parabolic PDEs that are far less degenerate
than equation (\ref{eq:1}). In particular, the results contained
in this paper can be easily extended to the case $\nu=0$, i.e. to
the evolutionary $p\,$-Poisson equation
\begin{equation}
u_{t}-\mathrm{div}\left(\vert Du\vert^{p-2}Du\right)=f\,\,\,\,\,\,\,\,\,\mathrm{in}\,\,\,\Omega_{T},\label{eq:p-Poisson}
\end{equation}
with $f\in L^{\frac{np\,+\,4}{np\,+\,4\,-n}}\left(0,T;W^{1,\frac{np\,+\,4}{np\,+\,4\,-n}}(\Omega)\right)$.
Therefore, our results permit to improve the existing literature,
already for equations of the form (\ref{eq:p-Poisson}), which exhibit
a milder degeneracy.
\end{singlespace}
\begin{singlespace}

\subsection{Motivation \label{subsec:Motivation}}
\end{singlespace}

\begin{singlespace}
\noindent $\hspace*{1em}$Before describing the structure of this
paper, we wish to motivate our study by stressing that, in the case
$n\leq3$, degenerate equations of the form (\ref{eq:1}) may arise
in \textit{gas filtration problems taking into account the initial
pressure gradient}. \\
$\hspace*{1em}$The existence of significant deviations from the linear
Darcy filtration law has been established for many systems consisting
of a fluid and a porous medium (e.g., the filtration of a gas in argillous
rocks). One of the manifestations of this nonlinearity is the existence
of a limiting (initial) pressure gradient, i.e. the minimum value
of the pressure gradient for which fluid motion occurs. In general,
fluid motion still takes place for subcritical values of the pressure
gradient, but very slowly; on reaching the limiting value of the pressure
gradient, there is a marked acceleration of the filtration. Therefore,
the limiting-gradient concept provides a good approximation for velocities
which are not too low.\\
$\hspace*{1em}$In accordance with some experimental results (see
\cite{Akh}), under certain physical conditions one can take the gas
filtration law in the very simple form
\[
\begin{cases}
\begin{array}{cc}
\mathbf{v}=-\,\frac{k}{\mu}\,D\mathrm{P}\left[1-\frac{G}{\left|D\mathrm{P}^{2}\right|}\right] & \,\,\,\mathrm{if}\,\,\left|D\mathrm{P}^{2}\right|\geq G,\\
\mathbf{v}=\mathbf{0}\,\,\,\,\,\,\,\,\,\,\,\,\,\,\,\,\,\,\,\,\,\,\,\,\,\,\,\,\,\,\,\,\,\,\,\,\,\,\,\,\,\,\,\,\, & \,\,\,\mathrm{if}\,\,\left|D\mathrm{P}^{2}\right|<G,
\end{array}\end{cases}
\]
where $\mathbf{v}=\mathbf{v}(x,t)$ is the filtration velocity, $k$
is the rock permeability, $\mu$ is the gas viscosity, $\mathrm{P}=\mathrm{P}(x,t)$
is the pressure and $G$ is a positive constant. Under this assumption
we obtain a particularly simple expression for the gas mass velocity
(flux) $\boldsymbol{j}$, which contains only the gradient of the
pressure squared, just as in the usual gas filtration problems:
\begin{equation}
\begin{cases}
\begin{array}{cc}
\boldsymbol{j}=\varrho\mathbf{v}=-\,\frac{k}{2\mu C}\left[D\mathrm{P}^{2}-G\,\frac{D\mathrm{P}^{2}}{\left|D\mathrm{P}^{2}\right|}\right] & \,\,\,\mathrm{if}\,\,\left|D\mathrm{P}^{2}\right|>G,\\
\boldsymbol{j}=\mathbf{0}\,\,\,\,\,\,\,\,\,\,\,\,\,\,\,\,\,\,\,\,\,\,\,\,\,\,\,\,\,\,\,\,\,\,\,\,\,\,\,\,\,\,\,\,\,\,\,\,\,\,\,\,\,\,\,\,\,\,\,\,\,\,\,\,\,\,\,\,\,\, & \,\,\,\mathrm{if}\,\,\left|D\mathrm{P}^{2}\right|\leq G,
\end{array}\end{cases}\label{eq:flux}
\end{equation}
where $\varrho$ is the gas density and $C$ is a positive constant.
Substituting expression (\ref{eq:flux}) into the gas mass-conservation
equation, we obtain the basic equation for the pressure:
\begin{equation}
\begin{cases}
\begin{array}{cc}
\frac{\partial\mathrm{P}}{\partial t}=\,\frac{k}{2m\mu}\,\mathrm{div}\left[D\mathrm{P}^{2}-G\,\frac{D\mathrm{P}^{2}}{\left|D\mathrm{P}^{2}\right|}\right] & \,\,\,\mathrm{if}\,\,\left|D\mathrm{P}^{2}\right|>G,\\
\frac{\partial\mathrm{P}}{\partial t}=0\,\,\,\,\,\,\,\,\,\,\,\,\,\,\,\,\,\,\,\,\,\,\,\,\,\,\,\,\,\,\,\,\,\,\,\,\,\,\,\,\,\,\,\,\,\,\,\,\,\,\,\,\,\,\,\,\,\,\,\, & \,\,\,\mathrm{if}\,\,\left|D\mathrm{P}^{2}\right|\leq G,
\end{array}\end{cases}\label{eq:pressure}
\end{equation}
where $m$ is a positive constant. Equation (\ref{eq:pressure}) implies,
first of all, that the steady gas motion is described by the same
relations as in the steady motion of an incompressible fluid if we
replace the pressure of the incompressible fluid with the square of
the gas pressure. Moreover, if the gas pressure differs very little
from some constant pressure $\mathrm{P}_{0}$, or if the gas pressure
differs considerably from a constant value only in regions where the
gas motion is nearly steady, then the equation for the gas filtration
in the region of motion can be ``linearized'' following L. S. Leibenson,
thus obtaining (see \cite{Akh} again):
\begin{equation}
\begin{cases}
\begin{array}{cc}
\frac{\partial\mathrm{P}^{2}}{\partial t}=\,\frac{k\,\mathrm{P}_{0}}{m\mu}\,\mathrm{div}\left[D\mathrm{P}^{2}-G\,\frac{D\mathrm{P}^{2}}{\left|D\mathrm{P}^{2}\right|}\right] & \,\,\,\mathrm{if}\,\,\left|D\mathrm{P}^{2}\right|>G,\\
\frac{\partial\mathrm{P}^{2}}{\partial t}=0\,\,\,\,\,\,\,\,\,\,\,\,\,\,\,\,\,\,\,\,\,\,\,\,\,\,\,\,\,\,\,\,\,\,\,\,\,\,\,\,\,\,\,\,\,\,\,\,\,\,\,\,\,\,\,\,\,\,\,\, & \,\,\,\mathrm{if}\,\,\left|D\mathrm{P}^{2}\right|\leq G.
\end{array}\end{cases}\label{eq:sqpress}
\end{equation}
$\hspace*{1em}$Now, setting $u=\mathrm{P}^{2}$ and performing an
appropriate scaling, equation (\ref{eq:sqpress}) turns into 
\[
\frac{\partial u}{\partial t}-\mathrm{div}\left[(\vert Du\vert-1)_{+}\,\,\frac{Du}{\vert Du\vert}\right]=0,
\]
which is nothing but equation (\ref{eq:1}) in the case $p=2$, $\nu=1$
and $f=0$. This is why (\ref{eq:1}) is sometimes called the \textit{Leibenson
equation} in the literature.\\
\\
$\hspace*{1em}$The paper is organized as follows. Section \ref{sec:prelim}
is devoted to the preliminaries: after a list of some classic notations
and some essentials estimates, we recall the basic properties of the
difference quotients of Sobolev functions. In Section \ref{sec:A priori},
we establish some \textit{a priori} estimates that will be needed
to demonstrate Theorem \ref{thm:main4}, whose proof is contained
in Section \ref{sec:main theo}. Using the existence of second spatial
derivatives, we then infer the existence of the weak time derivative
of the solutions. The corresponding arguments, which imply Theorem
\ref{thm:timeregularity}, are given in Section \ref{sec:timereg}.
\end{singlespace}
\begin{singlespace}

\section{Notations and preliminaries\label{sec:prelim}}
\end{singlespace}

\begin{singlespace}
\noindent $\hspace*{1em}$In this paper we shall denote by $C$ or
$c$ a general positive constant that may vary on different occasions.
Relevant dependencies on parameters and special constants will be
suitably emphasized using parentheses or subscripts. The norm we use
on $\mathbb{R}^{n}$ will be the standard Euclidean one and it will
be denoted by $\left|\,\cdot\,\right|$. In particular, for the vectors
$\xi,\eta\in\mathbb{R}^{n}$, we write $\langle\xi,\eta\rangle$ for
the usual inner product and $\left|\xi\right|:=\langle\xi,\xi\rangle^{\frac{1}{2}}$
for the corresponding Euclidean norm.\\
$\hspace*{1em}$For points in space-time, we will frequently use abbreviations
like $z=(x,t)$ or $z_{0}=(x_{0},t_{0})$, for spatial variables $x$,
$x_{0}\in\mathbb{R}^{n}$ and times $t$, $t_{0}\in\mathbb{R}$. We
also denote by $B(x_{0},\rho)=B_{\rho}(x_{0})=\left\{ x\in\mathbb{R}^{n}:\left|x-x_{0}\right|<\rho\right\} $
the open ball with radius $\rho>0$ and center $x_{0}\in\mathbb{R}^{n}$;
when not important, or clear from the context, we shall omit to denote
the center as follows: $B_{\rho}\equiv B(x_{0},\rho)$. Unless otherwise
stated, different balls in the same context will have the same center.
Moreover, we use the notation
\[
Q_{\rho}(z_{0}):=B_{\rho}(x_{0})\times(t_{0}-\rho^{2},t_{0}),\,\,\,\,\,z_{0}=(x_{0},t_{0})\in\mathbb{R}^{n}\times\mathbb{R},\,\,\rho>0,
\]
for the backward parabolic cylinder with vertex $(x_{0},t_{0})$ and
width $\rho$. We shall sometimes omit the dependence on the vertex
when all the cylinders occurring in a proof share the same vertex.
Finally, for a general cylinder $Q=A\times(t_{1},t_{2})$, where $A\subset\mathbb{R}^{n}$
and $t_{1}<t_{2}$, we denote by
\[
\partial_{\mathrm{par}}Q:=(\bar{A}\times\left\{ t_{1}\right\} )\cup(\partial A\times(t_{1},t_{2}))
\]
the usual \textit{parabolic boundary} of $Q$.\\
$\hspace*{1em}$We now recall some tools that will be useful to prove
our results. For the auxiliary function $H_{\lambda}:\mathbb{R}^{n}\rightarrow\mathbb{R}^{n}$
defined as
\[
H_{\lambda}(\xi):=\begin{cases}
\begin{array}{cc}
\left(\left|\xi\right|-\nu\right)_{+}^{\lambda}\frac{\xi}{\left|\xi\right|} & \,\,\,\,\,\,\,\,\,\,\,\,\,\,\,\,\,\,\,\mathrm{if}\,\,\,\xi\in\mathbb{R}^{n}\setminus\left\{ 0\right\} ,\\
0 & \mathrm{if}\,\,\,\xi=0,
\end{array}\end{cases}
\]
where $\lambda>0$ is a parameter, we record the following estimates,
which can be obtained by suitably modifying the proof of Lemma 4.1
in \cite{Br}.\\

\end{singlespace}
\begin{lem}
\begin{singlespace}
\noindent \label{lem:Brasco} If $2\leq p<\infty$, then for every
$\xi,\eta\in\mathbb{R}^{n}$ we get 
\[
\langle H_{p-1}(\xi)-H_{p-1}(\eta),\xi-\eta\rangle\,\geq\,\frac{4}{p^{2}}\left|H_{\frac{p}{2}}(\xi)-H_{\frac{p}{2}}(\eta)\right|^{2},
\]
\[
\left|H_{p-1}(\xi)-H_{p-1}(\eta)\right|\,\leq\,(p-1)\left(\left|H_{\frac{p}{2}}(\xi)\right|^{\frac{p-2}{p}}+\left|H_{\frac{p}{2}}(\eta)\right|^{\frac{p-2}{p}}\right)\left|H_{\frac{p}{2}}(\xi)-H_{\frac{p}{2}}(\eta)\right|.
\]
\end{singlespace}
\end{lem}

\begin{singlespace}
\noindent $\hspace*{1em}$In the following, we shall also use the
auxiliary function $V_{p}:\mathbb{R}^{n}\rightarrow\mathbb{R}^{n}$
defined as
\[
V_{p}(\xi):=\left|\xi\right|^{\frac{p-2}{2}}\xi,
\]
where $p\geq2$. For the above function, we recall the following estimates: 
\end{singlespace}
\begin{lem}
\begin{singlespace}
\noindent \label{lem:Lind} If $2\leq p<\infty$, then for every $\xi,\eta\in\mathbb{R}^{n}$
we get 
\[
\left|V_{p}(\xi)-V_{p}(\eta)\right|^{2}\,\leq\,\frac{p^{2}}{4}\,\langle\left|\xi\right|^{p-2}\xi-\left|\eta\right|^{p-2}\eta,\xi-\eta\rangle,
\]
\[
\left|\left|\xi\right|^{p-2}\xi-\left|\eta\right|^{p-2}\eta\right|\,\leq\,(p-1)\left(\left|\xi\right|^{\frac{p-2}{2}}+\left|\eta\right|^{\frac{p-2}{2}}\right)\left|V_{p}(\xi)-V_{p}(\eta)\right|.
\]
\end{singlespace}
\end{lem}

\begin{singlespace}
\noindent We refer to \cite[Chapter 12]{Lind} for a proof of these
fundamental inequalities. \\
$\hspace*{1em}$The next lemma has an important application in the
so called hole-filling method and its proof can be found, for example,
in \cite[Lemma 6.1]{Giu}.
\end{singlespace}
\begin{lem}
\begin{singlespace}
\noindent \label{lem:giusti} Assume that $\Psi:[r_{0},r_{1}]\rightarrow[0,\infty)$
is a bounded function which satisfies 
\[
\Psi(s)\,\leq\,\theta\,\Psi(t)\,+\,\frac{A}{(t-s)^{\alpha}}\,+\,\frac{B}{(t-s)^{\beta}}\,+\,C
\]
for all $r_{0}\leq s<t\leq r_{1}$ and fixed non-negative constants
$A$, $B$, $C$, $\alpha\geq\beta>0$ and $\theta\in(0,1)$. Then
\[
\Psi(r_{0})\,\leq\,c\,\left(\frac{A}{(r_{1}-r_{0})^{\alpha}}\,+\,\frac{B}{(r_{1}-r_{0})^{\beta}}\,+\,C\right),
\]
where $c\equiv c(\alpha,\theta)>0$.
\end{singlespace}
\end{lem}

\begin{singlespace}
\noindent $\hspace*{1em}$For further needs, we also record the following
interpolation inequality, whose proof can be found in \cite[Proposition 3.1]{DiBene}.
\end{singlespace}
\begin{lem}
\begin{singlespace}
\noindent \label{lem:inter} Assume that the function $v:Q_{r}(z_{0})\cup\partial_{\mathrm{par}}Q_{r}(z_{0})\rightarrow\mathbb{R}$
satisfies
\[
v\in L^{\infty}\left(t_{0}-r^{2},t_{0};L^{q}\left(B_{r}(x_{0})\right)\right)\cap L^{p}\left(t_{0}-r^{2},t_{0};W_{0}^{1,p}\left(B_{r}(x_{0})\right)\right)
\]

\noindent for some exponents $1\leq p$,$q<\infty$. Then the following
estimate
\[
\int_{Q_{r}(z_{0})}\left|v\right|^{p\,+\,pq/n}dz\,\leq\,c\left(\sup_{s\in(t_{0}-r^{2},t_{0})}\int_{B_{r}(x_{0})}\left|v(x,s)\right|^{q}dx\right)^{p/n}\int_{Q_{r}(z_{0})}\left|Dv\right|^{p}dz
\]
holds true for a positive constant $c$ depending at most on $n$,
$p$ and $q$.
\end{singlespace}
\end{lem}

\begin{singlespace}
\noindent $\hspace*{1em}$We conclude by recalling the following
\end{singlespace}
\begin{defn}
\begin{singlespace}
\noindent A function $u\in C^{0}\left((0,T);L^{2}\left(\Omega\right)\right)\cap L^{p}\left(0,T;W^{1,p}\left(\Omega\right)\right)$
is a \textit{weak solution} of equation (\ref{eq:1}) if and only
if for any test function $\varphi\in C_{0}^{\infty}(\Omega_{T})$
the following integral identity holds:
\begin{equation}
\int_{\Omega_{T}}\left(u\cdot\partial_{t}\varphi-\langle H_{p-1}(Du),D\varphi\rangle\right)\,dz\,=\,-\int_{\Omega_{T}}f\varphi\,dz.\label{eq:weaksol}
\end{equation}
\end{singlespace}
\end{defn}

\begin{singlespace}

\subsection{Difference quotients\label{subsec:DiffOpe}}
\end{singlespace}

\begin{singlespace}
\noindent $\hspace*{1em}$We recall here the definition and some elementary
properties of the difference quotients that will be useful in the
following (see, for example, \cite{Giu}).\\

\end{singlespace}
\begin{defn}
\begin{singlespace}
\noindent For every vector-valued function $F:\mathbb{R}^{n}\rightarrow\mathbb{R}^{N}$
the \textit{finite difference operator }in the direction $x_{s}$
is defined by
\[
\tau_{s,h}F(x)=F(x+he_{s})-F(x),
\]
where $h\in\mathbb{R}$, $e_{s}$ is the unit vector in the direction
$x_{s}$ and $s\in\left\{ 1,\ldots,n\right\} $. \\
$\hspace*{1em}$The \textit{difference quotient} of $F$ with respect
to $x_{s}$ is defined for $h\in\mathbb{R}\setminus\left\{ 0\right\} $
as 
\[
\Delta_{s,h}F(x)=\frac{\tau_{s,h}F(x)}{h}.
\]
\end{singlespace}
\end{defn}

\begin{singlespace}
\noindent $\hspace*{1em}$When no confusion can arise, we shall omit
the index $s$ and simply write $\tau_{h}$ or $\Delta_{h}$ instead
of $\tau_{s,h}$ or $\Delta_{s,h}$, respectively. 
\end{singlespace}
\begin{prop}
\begin{singlespace}
\noindent Let $F$ be a function such that $F\in W^{1,p}\left(\Omega\right)$,
with $p\geq1$, and let us consider the set
\[
\Omega_{\left|h\right|}:=\left\{ x\in\Omega:\mathrm{dist}\left(x,\partial\Omega\right)>\left|h\right|\right\} .
\]
Then:\\
\\
$\mathrm{(}a\mathrm{)}$ $\Delta_{h}F\in W^{1,p}\left(\Omega_{\left|h\right|}\right)$
and $D_{i}(\Delta_{h}F)=\Delta_{h}(D_{i}F)$ for every $\,i\in\left\{ 1,\ldots,n\right\} $.\\

\noindent $\mathrm{(}b\mathrm{)}$ If at least one of the functions
$F$ or $G$ has support contained in $\Omega_{\left|h\right|}$,
then
\[
\int_{\Omega}F\,\Delta_{h}G\,dx\,=\,-\int_{\Omega}G\,\Delta_{-h}F\,dx.
\]
$\mathrm{(}c\mathrm{)}$ We have 
\[
\Delta_{h}(FG)(x)=F(x+he_{s})\Delta_{h}G(x)\,+\,G(x)\Delta_{h}F(x).
\]
\end{singlespace}
\end{prop}

\begin{singlespace}
\noindent $\hspace*{1em}$The next result about the finite difference
operator is a kind of integral version of Lagrange Theorem and can
be obtained by combining Lemma 8.1 in \cite{Giu} with the theorem
on page 3 of \cite{Maz'ya}. 
\end{singlespace}
\begin{lem}
\begin{singlespace}
\noindent \label{lem:Giusti1} If $0<\rho<R$, $\left|h\right|<\frac{R-\rho}{2}$,
$1<q<+\infty$ and $F\in L_{loc}^{1}\left(B_{R},\mathbb{R}^{N}\right)$
is such that $DF\in L^{q}\left(B_{R},\mathbb{R}^{N\times n}\right)$,
then
\[
\int_{B_{\rho}}\left|\tau_{h}F(x)\right|^{q}dx\,\leq\,c^{q}(n)\left|h\right|^{q}\int_{B_{R}}\left|DF(x)\right|^{q}dx.
\]
Moreover, if $F\in L^{q}\left(B_{R},\mathbb{R}^{N}\right)$ then we
have 
\[
\int_{B_{\rho}}\left|F(x+he_{s})\right|^{q}dx\,\leq\,\int_{B_{R}}\left|F(x)\right|^{q}dx.
\]
\end{singlespace}
\end{lem}

\begin{singlespace}
\noindent $\hspace*{1em}$Finally, we recall the following fundamental
result, whose proof can be found in \cite[Lemma 8.2]{Giu}:
\end{singlespace}
\begin{lem}
\begin{singlespace}
\noindent \label{lem:RappIncre} Let $F:\mathbb{R}^{n}\rightarrow\mathbb{R}^{N}$,
$F\in L^{q}\left(B_{R},\mathbb{R}^{N}\right)$ with $1<q<+\infty$.
Suppose that there exist $\rho\in(0,R)$ and a constant $M>0$ such
that 
\[
\sum_{s=1}^{n}\int_{B_{\rho}}\left|\tau_{s,h}F(x)\right|^{q}dx\,\leq\,M^{q}\left|h\right|^{q}
\]
for every $h$ with $\left|h\right|<\frac{R-\rho}{2}$. Then $F\in W^{1,q}\left(B_{\rho},\mathbb{R}^{N}\right)$.
Moreover 
\[
\Vert DF\Vert_{L^{q}\left(B_{\rho}\right)}\leq M
\]
and
\[
\Delta_{s,h}F\rightarrow D_{s}F\,\,\,\,\,\,\,\,\,\,in\,\,L_{loc}^{q}\left(B_{R}\right),\,\,as\,\,h\rightarrow0,
\]
for each $s\in\left\{ 1,\ldots,n\right\} $.
\end{singlespace}
\end{lem}

\begin{singlespace}

\section{\textit{A priori} estimates\label{sec:A priori}}
\end{singlespace}

\begin{singlespace}
\noindent $\hspace*{1em}$In this section, we shall derive some \textit{a
priori} estimates for the solutions $u_{\varepsilon}$ of equation
(\ref{eq:regequ}) below, which is much less degenerate than equation
(\ref{eq:1}). Such estimates will play a key role in the proof of
Theorem \ref{thm:main4} (see Section \ref{sec:main theo}).\\
$\hspace*{1em}$Let $n\geq2$, $p\geq2$, $\frac{np\,+\,4}{np\,+\,4\,-n}\leq\vartheta<\infty$
and $f\in L^{\vartheta}\left(0,T;W^{1,\vartheta}(\Omega)\right)$.
For $\varepsilon\in[0,1]$ and a couple of standard, non-negative,
radially symmetric mollifiers $\phi_{1}\in C_{0}^{\infty}(B_{1}(0))$
and $\phi_{2}\in C_{0}^{\infty}((-1,1))$ we write
\begin{equation}
f_{\varepsilon}(x,t):=\int_{-1}^{1}\int_{B_{1}(0)}f(x-\varepsilon y,t-\varepsilon s)\,\phi_{1}(y)\,\phi_{2}(s)\,dy\,ds,\label{eq:molli}
\end{equation}
where $f$ is meant to be extended by zero outside $\Omega_{T}$.
Let us observe that $f_{0}=f$ and $f_{\varepsilon}\in C^{\infty}(\Omega_{T})$
for every $\varepsilon\in(0,1]$.\\
Now we consider a domain in space-time denoted by $\Omega'_{1,2}:=\Omega'\times(T_{1},T_{2})$,
where $\Omega'\subseteq\Omega$ is a bounded domain with smooth boundary
and $(T_{1},T_{2})\Subset(0,T)$. For our purposes, in the following
we will need the definitions below.
\end{singlespace}
\begin{defn}
\begin{singlespace}
\noindent Let $\varepsilon\in(0,1]$. A function $u_{\varepsilon}\in C^{0}\left((T_{1},T_{2});L^{2}(\Omega')\right)\cap L^{p}\left(T_{1},T_{2};W^{1,p}(\Omega')\right)$
is a \textit{weak solution} of the equation 
\begin{equation}
\partial_{t}u_{\varepsilon}-\mathrm{div}\left(H_{p-1}(Du_{\varepsilon})+\varepsilon\,\vert Du_{\varepsilon}\vert^{p-2}Du_{\varepsilon}\right)=f_{\varepsilon}\,\,\,\,\,\,\mathrm{in}\,\,\,\Omega'_{1,2}\label{eq:regequ}
\end{equation}
if and only if for any test function $\varphi\in C_{0}^{\infty}(\Omega'_{1,2})$
the following integral identity holds:
\begin{equation}
\int_{\Omega'_{1,2}}\left(u_{\varepsilon}\cdot\partial_{t}\varphi-\langle H_{p-1}(Du_{\varepsilon})+\varepsilon\,\vert Du_{\varepsilon}\vert^{p-2}Du_{\varepsilon},D\varphi\rangle\right)\,dz\,=\,-\int_{\Omega'_{1,2}}f_{\varepsilon}\,\varphi\,dz.\label{eq:weaksol2}
\end{equation}
\end{singlespace}
\end{defn}

\begin{singlespace}
\noindent 

\end{singlespace}\begin{defn}
\begin{singlespace}
\noindent \label{def:L2-traces} Let $\varepsilon\in(0,1]$ and $g\in C^{0}\left([T_{1},T_{2}];L^{2}(\Omega')\right)\cap L^{p}\left(T_{1},T_{2};W^{1,p}(\Omega')\right)$.
In this framework, we identify a function 
\[
u_{\varepsilon}\in C^{0}\left([T_{1},T_{2}];L^{2}(\Omega')\right)\cap L^{p}\left(T_{1},T_{2};W^{1,p}(\Omega')\right)
\]
as a \textit{weak solution of the Cauchy-Dirichlet problem\medskip{}
}
\begin{equation}
\begin{cases}
\begin{array}{cc}
\partial_{t}u_{\varepsilon}-\mathrm{div}\left(H_{p-1}(Du_{\varepsilon})+\varepsilon\,\vert Du_{\varepsilon}\vert^{p-2}Du_{\varepsilon}\right)=f_{\varepsilon} & \,\,\,\mathrm{in}\,\,\,\Omega'_{1,2},\\
u_{\varepsilon}=g & \,\,\,\,\,\,\,\,\,\,\,\,\,\,\mathrm{on}\,\,\,\partial_{\mathrm{par}}\Omega'_{1,2},
\end{array}\end{cases}\label{eq:CAUCHYDIR}
\end{equation}
if and only if (\ref{eq:weaksol2}) holds and moreover, $u_{\varepsilon}\in g+L^{p}\left(T_{1},T_{2};W_{0}^{1,p}(\Omega')\right)$
and $u_{\varepsilon}(\cdot,T_{1})=g(\cdot,T_{1})$ in the $L^{2}$-sense,
that is 
\begin{equation}
\underset{t\,\rightarrow\,(T_{1})^{+}}{\lim}\,\Vert u_{\varepsilon}(\cdot,t)-g(\cdot,T_{1})\Vert_{L^{2}(\Omega')}\,=\,0.\label{eq:L2sense}
\end{equation}
Therefore, the initial condition $u_{\varepsilon}=g$ on $\Omega'\times\{T_{1}\}$
has to be understood in the usual $L^{2}$-sense (\ref{eq:L2sense}),
while the condition $u_{\varepsilon}=g$ on the lateral boundary $\partial\Omega'\times(T_{1},T_{2})$
has to be meant in the sense of traces, i.e. $(u_{\varepsilon}-g)\,(\cdot,t)\in W_{0}^{1,p}\left(\Omega'\right)$
for almost every $t\in(T_{1},T_{2})$.
\end{singlespace}
\end{defn}

\begin{singlespace}
\noindent \begin{brem} The regularized parabolic equation (\ref{eq:regequ})
fulfills standard $p\,$-growth conditions. The advantage of considering
the associated Cauchy-Dirichlet problem (\ref{eq:CAUCHYDIR}) arises
from the fact that the existence of a unique solution $u_{\varepsilon}\in C^{0}\left([T_{1},T_{2}];L^{2}(\Omega')\right)\cap L^{p}\left(T_{1},T_{2};W^{1,p}(\Omega')\right)$
satisfying the requirements of Definition \ref{def:L2-traces} can
be ensured by the classic existence theory for parabolic equations (see \cite[Chapter 2, Theorem 1.2 and Remark 1.2]{Lions}).\end{brem}\medskip{}

\noindent $\hspace*{1em}$Now we shall use the well-known difference
quotients method in the spatial directions together with the properties
of the functions $H_{\lambda}$ and $V_{p}$ to establish the following
result.
\end{singlespace}
\begin{prop}
\begin{singlespace}
\noindent \label{prop:apriori1} Let $\varepsilon\in(0,1]$, $n\geq2$,
$p\geq2$, $\frac{np\,+\,4}{np\,+\,4\,-n}\leq\vartheta<\infty$ and
$f\in L^{\vartheta}\left(0,T;W^{1,\vartheta}(\Omega)\right)$. Moreover,
assume that 
\[
u_{\varepsilon}\in C^{0}\left((T_{1},T_{2});L^{2}\left(\Omega'\right)\right)\cap L^{p}\left(T_{1},T_{2};W^{1,p}\left(\Omega'\right)\right)
\]
is a weak solution of equation $\mathrm{(\ref{eq:regequ})}$. Then
the solution satisfies 
\begin{equation}
H_{\frac{p}{2}}(Du_{\varepsilon})\,\in\,L_{loc}^{2}(T_{1},T_{2};W_{loc}^{1,2}\left(\Omega',\mathbb{R}^{n}\right))\,\,\,\,\,\,\,\,\,\,\,\,\,\,and\,\,\,\,\,\,\,\,\,\,\,\,\,\,Du_{\varepsilon}\,\in\,L_{loc}^{\infty}(T_{1},T_{2};L_{loc}^{2}\left(\Omega',\mathbb{R}^{n}\right)).\label{eq:apriori1}
\end{equation}
Furthermore, the following estimate\begin{equation}\label{eq:Caccio1}
\begin{split}
\sup_{t\in(t_{0}-(\rho/2)^{2},t_{0})}&\Vert Du_{\varepsilon}(\cdot,t)\Vert_{L^{2}(B_{\rho/2}(x_{0}))}^{2}\,+\int_{Q_{\rho/2}(z_{0})}\left|DH_{\frac{p}{2}}(Du_{\varepsilon})\right|^{2}dz\\
\,\,\,\,\,\,\,\,\,\,\,\,\,\,&\leq\,c\left[\Vert Df_{\varepsilon}\Vert_{L^{p'}(Q_{R_{0}})}\left(\int_{Q_{R_{0}}}\left|Du_{\varepsilon}\right|^{p}dz\right)^{\frac{1}{p}}+\,\rho^{-2}\int_{Q_{R_{0}}}\left(\left|Du_{\varepsilon}\right|^{p}+1\right)\,dz\right]
\end{split}
\end{equation}holds true for any parabolic cylinder $Q_{\rho}(z_{0})\subset Q_{R_{0}}(z_{0})\Subset\Omega'_{1,2}$
and a positive constant $c$ depending only on $n$ and $p$.
\end{singlespace}
\end{prop}

\begin{singlespace}
\noindent \begin{proof}[\bfseries{Proof}] By a slight abuse of notation,
for $w\in L_{loc}^{1}\left(\Omega'_{1,2},\mathbb{R}^{N}\right)$ and
$s\in\left\{ 1,\ldots,n\right\} $, $h\neq0$, we set (when $x+he_{s}\in\Omega'$)
\[
\tau_{h}w(x,t)\equiv\tau_{s,h}w(x,t):=w(x+he_{s},t)-w(x,t),
\]
\[
\Delta_{h}w(x,t)\equiv\Delta_{s,h}w(x,t):=\frac{w(x+he_{s},t)-w(x,t)}{h},
\]
\\
where $e_{s}$ is the unit vector in the direction $x_{s}$. \\
$\hspace*{1em}$Since $u_{\varepsilon}$ is a weak solution of equation
(\ref{eq:regequ}), we have\\
\[
\int_{\Omega'_{1,2}}\left(u_{\varepsilon}\cdot\partial_{t}\varphi-\langle H_{p-1}(Du_{\varepsilon})+\varepsilon\,\vert Du_{\varepsilon}\vert^{p-2}Du_{\varepsilon},D\varphi\rangle\right)\,dz\,=\,-\int_{\Omega'_{1,2}}f_{\varepsilon}\,\varphi\,dz,
\]
\\
for every test function $\varphi\in C_{0}^{\infty}(\Omega'_{1,2})$.
Replacing $\varphi$ with $\tau_{-h}\varphi$, where $0<\left|h\right|<\mathrm{dist}(\mathrm{supp}\,\varphi,\partial\Omega'_{1,2})$,
by virtue of the properties of the finite difference operator, we
get\\
\[
\int_{\Omega'_{1,2}}\left(\tau_{h}u_{\varepsilon}\cdot\partial_{t}\varphi-\langle\tau_{h}H_{p-1}(Du_{\varepsilon})+\varepsilon\,\tau_{h}\left[\vert Du_{\varepsilon}\vert^{p-2}Du_{\varepsilon}\right],D\varphi\rangle\right)\,dz\,=\,-\int_{\Omega'_{1,2}}\tau_{h}f_{\varepsilon}\cdot\varphi\,dz.
\]
\\
We now replace $\varphi$ by $\varphi_{\sigma}\equiv\phi_{\sigma}\ast\varphi$
in the previous equation, where $\left\{ \phi_{\sigma}\right\} $,
$\sigma>0$, denotes the family of standard, non-negative, radially
symmetric mollifiers in $\mathbb{R}^{n+1}$. This yields, for $0<\sigma\ll1$\\
\[
\int_{\Omega'_{1,2}}\left((\tau_{h}u_{\varepsilon})_{\sigma}\cdot\partial_{t}\varphi-\langle(\tau_{h}H_{p-1}(Du_{\varepsilon}))_{\sigma}+\varepsilon\left(\tau_{h}\left[\vert Du_{\varepsilon}\vert^{p-2}Du_{\varepsilon}\right]\right)_{\sigma},D\varphi\rangle\right)\,dz\,=\,-\int_{\Omega'_{1,2}}(\tau_{h}f_{\varepsilon})_{\sigma}\cdot\varphi\,dz.
\]
\\
Now, in the last equation we choose the test function $\varphi\equiv\Phi(\tau_{h}u_{\varepsilon})_{\sigma}$,
where $\Phi\in C_{0}^{\infty}(\Omega'_{1,2})$ is a smooth function
which will be specified later. After an integration by parts and then
letting $\sigma\rightarrow0$, we obtain\\
 \begin{equation}\label{eq:rif1}
\begin{split}
&-\,\frac{1}{2}\int_{\Omega'_{1,2}}\left|\tau_{h}u_{\varepsilon}\right|^{2}\partial_{t}\Phi\,\,dz\,+\int_{\Omega'_{1,2}}\Phi\,\langle\tau_{h}H_{p-1}(Du_{\varepsilon})+\,\varepsilon\,\tau_{h}\left[\vert Du_{\varepsilon}\vert^{p-2}Du_{\varepsilon}\right],D\tau_{h}u_{\varepsilon}\rangle\,dz\\
&=\,-\int_{\Omega'_{1,2}}\langle\tau_{h}H_{p-1}(Du_{\varepsilon}),D\Phi\rangle\tau_{h}u_{\varepsilon}\,\,dz\,-\varepsilon\int_{\Omega'_{1,2}}\langle\tau_{h}\left[\vert Du_{\varepsilon}\vert^{p-2}Du_{\varepsilon}\right],D\Phi\rangle\tau_{h}u_{\varepsilon}\,\,dz\\
&\,\,\,\,\,\,\,+\int_{\Omega'_{1,2}}\tau_{h}f_{\varepsilon}\cdot\Phi\cdot\tau_{h}u_{\varepsilon}\,\,dz,
\end{split}
\end{equation}for any smooth function $\Phi\in C_{0}^{\infty}(\Omega'_{1,2})$,
provided that $\left|h\right|$ is small enough. Note that an approximation
argument yields the same identity for any $\Phi\in W^{1,\infty}(\Omega'_{1,2})$
with compact support in $\Omega'_{1,2}$ and any sufficiently small
$h\in\mathbb{R}\setminus\left\{ 0\right\} $. In what follows, we
will denote by $c_{k}$ and $c$ some positive constants which do
not depend either on $h$ or $\varepsilon$.\\
$\hspace*{1em}$Now, let us consider a parabolic cylinder $Q_{\rho}(z_{0})\subset Q_{R_{0}}(z_{0})\Subset\Omega'_{1,2}$.
For a fixed time $t_{1}\in(t_{0}-\rho^{2},t_{0})$ and $\delta\in(0,t_{0}-t_{1})$,
we choose $\Phi(x,t)=\widetilde{\chi}(t)\chi(t)\eta^{2}(x)$ with
$\chi\in W^{1,\infty}\left((T_{1},T_{2}),\left[0,1\right]\right)$,
$\chi\equiv0$ on $(T_{1},t_{0}-\rho^{2})$ and $\partial_{t}\chi\geq0$,
$\eta\in C_{0}^{\infty}\left(B_{\rho}(x_{0}),\left[0,1\right]\right)$,
and with the Lipschitz continuous function $\widetilde{\chi}:(T_{1},T_{2})\rightarrow\mathbb{R}$
defined by
\[
\widetilde{\chi}(t)=\begin{cases}
\begin{array}{c}
1\\
\mathrm{affine}\\
0
\end{array} & \begin{array}{c}
\mathrm{if}\,\,t\leq t_{1},\quad\quad\quad\;\,\,\\
\mathrm{if}\,\,t_{1}<t<t_{1}+\delta,\\
\mathrm{if}\,\,t\geq t_{1}+\delta.\quad\quad
\end{array}\end{cases}
\]

\noindent With such a choice of $\Phi$, equation \eqref{eq:rif1}
becomes\\
\begin{align*}
&-\,\frac{1}{2}\int_{\Omega'_{1,2}}\left|\tau_{h}u_{\varepsilon}\right|^{2}\eta^{2}(x)\chi(t)\partial_{t}\widetilde{\chi}(t)\,\,dz\,-\,\frac{1}{2}\int_{\Omega'_{1,2}}\left|\tau_{h}u_{\varepsilon}\right|^{2}\eta^{2}(x)\widetilde{\chi}(t)\partial_{t}\chi(t)\,\,dz\\
&+\int_{\Omega'_{1,2}}\widetilde{\chi}(t)\chi(t)\eta^{2}(x)\langle\tau_{h}H_{p-1}(Du_{\varepsilon}),D\tau_{h}u_{\varepsilon}\rangle\,dz\,+\varepsilon\int_{\Omega'_{1,2}}\widetilde{\chi}(t)\chi(t)\eta^{2}(x)\langle\tau_{h}[\vert Du_{\varepsilon}\vert^{p-2}Du_{\varepsilon}],D\tau_{h}u_{\varepsilon}\rangle\,dz\\
&\,\,\,\,=-2\int_{\Omega'_{1,2}}\widetilde{\chi}(t)\chi(t)\eta(x)\langle\tau_{h}H_{p-1}(Du_{\varepsilon}),D\eta\rangle\tau_{h}u_{\varepsilon}\,\,dz\\
&\,\,\,\,\,\,\,\,\,\,-2\,\varepsilon\int_{\Omega'_{1,2}}\widetilde{\chi}(t)\chi(t)\eta(x)\langle\tau_{h}[\vert Du_{\varepsilon}\vert^{p-2}Du_{\varepsilon}],D\eta\rangle\tau_{h}u_{\varepsilon}\,\,dz\,+\int_{\Omega'_{1,2}}(\tau_{h}f_{\varepsilon})(\tau_{h}u_{\varepsilon})\widetilde{\chi}(t)\chi(t)\eta^{2}(x)\,dz.
\end{align*}Setting $Q^{t_{1}}:=B_{\rho}(x_{0})\times(t_{0}-\rho^{2},t_{1})$
and letting $\delta\rightarrow0$ in the previous equality, for every
$t_{1}\in(t_{0}-\rho^{2},t_{0})$ we get\begin{equation}\label{eq:rif2}
\begin{split}
&\frac{1}{2}\int_{B_{\rho}(x_{0})}\chi(t_{1})\eta^{2}(x)\left|\tau_{h}u_{\varepsilon}(x,t_{1})\right|^{2}dx\,+\int_{Q^{t_{1}}}\chi(t)\eta^{2}(x)\langle\tau_{h}H_{p-1}(Du_{\varepsilon}),D\tau_{h}u_{\varepsilon}\rangle\,dz\\
&+\,\varepsilon\int_{Q^{t_{1}}}\chi(t)\eta^{2}(x)\langle\tau_{h}[\vert Du_{\varepsilon}\vert^{p-2}Du_{\varepsilon}],D\tau_{h}u_{\varepsilon}\rangle\,dz\\
&=-2\int_{Q^{t_{1}}}\chi(t)\eta(x)\langle\tau_{h}H_{p-1}(Du_{\varepsilon}),D\eta\rangle\tau_{h}u_{\varepsilon}\,\,dz\,+\int_{Q^{t_{1}}}(\tau_{h}f_{\varepsilon})(\tau_{h}u_{\varepsilon})\chi(t)\eta^{2}(x)\,dz\\
&\,\,\,\,\,\,-2\,\varepsilon\int_{Q^{t_{1}}}\chi(t)\eta(x)\langle\tau_{h}[\vert Du_{\varepsilon}\vert^{p-2}Du_{\varepsilon}],D\eta\rangle\tau_{h}u_{\varepsilon}\,\,dz\,+\,\frac{1}{2}\int_{Q^{t_{1}}}(\partial_{t}\chi)\eta^{2}(x)\left|\tau_{h}u_{\varepsilon}\right|^{2}dz\\
&=:I_{1}+I_{2}+I_{3}+I_{4},
\end{split}
\end{equation}\\
where we have used that $\partial_{t}\widetilde{\chi}$ converges
to a Dirac delta distribution as $\delta\rightarrow0$, together with
the $L^{2}(\Omega')-$valued continuity of $u_{\varepsilon}$. In
the following, we estimate $I_{1}$, $I_{2}$ and $I_{3}$ separately.
Let us first consider $I_{1}$. Using Lemma \ref{lem:Brasco} together
with Young's inequality, we obtain\begin{equation}\label{eq:rif3}
\begin{split}
\left|I_{1}\right|\leq\,&\,\,\frac{2(p-1)^{2}}{\beta}\int_{Q^{t_{1}}}\chi(t)\left|D\eta\right|^{2}\left(\left|H_{\frac{p}{2}}(Du_{\varepsilon}(x+he_{s},t))\right|^{\frac{p-2}{p}}+\left|H_{\frac{p}{2}}(Du_{\varepsilon})\right|^{\frac{p-2}{p}}\right)^{2}\left|\tau_{h}u_{\varepsilon}\right|^{2}dz\\
&+\,\frac{\beta}{2}\int_{Q^{t_{1}}}\chi(t)\eta^{2}(x)\left|\tau_{h}H_{\frac{p}{2}}(Du_{\varepsilon})\right|^{2}dz,
\end{split}
\end{equation}where $\beta>0$ will be chosen later. As for $I_{3}$, by Lemma \ref{lem:Lind}
we similarly have\\
\begin{equation}\label{eq:rif4}
\begin{split}
\left|I_{3}\right|\leq\,&\,\,\frac{2\varepsilon\,(p-1)^{2}}{\beta}\int_{Q^{t_{1}}}\chi(t)\left|D\eta\right|^{2}\left(\left|Du_{\varepsilon}(x+he_{s},t)\right|^{\frac{p-2}{2}}+\left|Du_{\varepsilon}\right|^{\frac{p-2}{2}}\right)^{2}\left|\tau_{h}u_{\varepsilon}\right|^{2}dz\\
&+\,\frac{\varepsilon\beta}{2}\int_{Q^{t_{1}}}\chi(t)\eta^{2}(x)\left|\tau_{h}V_{p}(Du_{\varepsilon})\right|^{2}dz.
\end{split}
\end{equation}\\
We now turn our attention to $I_{2}$. Taking advantage of the properties
of $f_{\varepsilon}$, $u_{\varepsilon}$, $\chi$ and $\eta$, and
using Hölder's inequality with exponents $p$ and $p'$ together with
Lemma \ref{lem:Giusti1}, we obtain\\
\begin{equation}
\left|I_{2}\right|\,\leq\,c_{1}(n)\,\left|h\right|^{2}\Vert Df_{\varepsilon}\Vert_{L^{p'}(Q_{R_{0}})}\,\Vert Du_{\varepsilon}\Vert_{L^{p}(Q_{R_{0}})},\label{eq:rif5}
\end{equation}
provided that $\left|h\right|$ is suitably small. Now, by virtue
of Lemma \ref{lem:Brasco} we also have\\
\begin{equation}
\frac{4}{p^{2}}\int_{Q^{t_{1}}}\chi(t)\eta^{2}(x)\left|\tau_{h}H_{\frac{p}{2}}(Du_{\varepsilon})\right|^{2}dz\,\leq\,\int_{Q^{t_{1}}}\chi(t)\eta^{2}(x)\langle\tau_{h}H_{p-1}(Du_{\varepsilon}),D\tau_{h}u_{\varepsilon}\rangle\,dz.\label{eq:rif6}
\end{equation}
\\
Similarly, by Lemma \ref{lem:Lind} we obtain\\
\begin{equation}
\frac{4\,\varepsilon}{p^{2}}\int_{Q^{t_{1}}}\chi(t)\eta^{2}(x)\left|\tau_{h}V_{p}(Du_{\varepsilon})\right|^{2}dz\,\leq\,\varepsilon\int_{Q^{t_{1}}}\chi(t)\eta^{2}(x)\langle\tau_{h}[\vert Du_{\varepsilon}\vert^{p-2}Du_{\varepsilon}],D\tau_{h}u_{\varepsilon}\rangle\,dz.\label{eq:rif7}
\end{equation}
\\
Collecting estimates \eqref{eq:rif2}, \eqref{eq:rif3}, \eqref{eq:rif4},
(\ref{eq:rif5}), (\ref{eq:rif6}) and (\ref{eq:rif7}), and choosing
$\beta=4/p^{2}$, we arrive at\\
\begin{align*}
&\int_{B_{\rho}(x_{0})}\chi(t_{1})\eta^{2}(x)\left|\tau_{h}u_{\varepsilon}(x,t_{1})\right|^{2}dx\,+\int_{Q^{t_{1}}}\chi(t)\eta^{2}(x)\left|\tau_{h}H_{\frac{p}{2}}(Du_{\varepsilon})\right|^{2}dz\\
&+\,\,\varepsilon\int_{Q^{t_{1}}}\chi(t)\eta^{2}(x)\left|\tau_{h}V_{p}(Du_{\varepsilon})\right|^{2}dz\\
&\leq\,c_{2}(n,p)\int_{Q^{t_{1}}}\left[\chi\left|D\eta\right|^{2}\left(\left|H_{\frac{p}{2}}(Du_{\varepsilon}(x+he_{s},t))\right|^{\frac{p-2}{p}}+\left|H_{\frac{p}{2}}(Du_{\varepsilon})\right|^{\frac{p-2}{p}}\right)^{2}+(\partial_{t}\chi)\eta^{2}\right]\left|\tau_{h}u_{\varepsilon}\right|^{2}dz\\
&\,\,\,\,\,\,\,+\,c_{2}(n,p)\int_{Q^{t_{1}}}\left[\varepsilon\,\chi\left|D\eta\right|^{2}\left(\left|Du_{\varepsilon}(x+he_{s},t)\right|^{\frac{p-2}{2}}+\left|Du_{\varepsilon}\right|^{\frac{p-2}{2}}\right)^{2}+(\partial_{t}\chi)\eta^{2}\right]\left|\tau_{h}u_{\varepsilon}\right|^{2}dz\\
&\,\,\,\,\,\,\,+\,\,c_{2}(n,p)\left|h\right|^{2}\Vert Df_{\varepsilon}\Vert_{L^{p'}(Q_{R_{0}})}\,\Vert Du_{\varepsilon}\Vert_{L^{p}(Q_{R_{0}})},
\end{align*}which holds for every $t_{1}\in(t_{0}-\rho^{2},t_{0})$ and every
sufficiently small $h\in\mathbb{R}\setminus\left\{ 0\right\} $. Recalling
that $0<\varepsilon\leq1$ and using the definition of the function
$H_{\frac{p}{2}}$, from the above estimate we get\\
 \begin{align*}
&\int_{B_{\rho}(x_{0})}\chi(t_{1})\eta^{2}(x)\left|\tau_{h}u_{\varepsilon}(x,t_{1})\right|^{2}dx\,+\int_{Q^{t_{1}}}\chi(t)\eta^{2}(x)\left|\tau_{h}H_{\frac{p}{2}}(Du_{\varepsilon})\right|^{2}dz\\
&\leq\,c_{3}(n,p)\int_{Q^{t_{1}}}\left[\chi\left|D\eta\right|^{2}\left(\left|Du_{\varepsilon}(x+he_{s},t)\right|^{\frac{p-2}{2}}+\left|Du_{\varepsilon}\right|^{\frac{p-2}{2}}\right)^{2}+(\partial_{t}\chi)\eta^{2}\right]\left|\tau_{h}u_{\varepsilon}\right|^{2}dz\\
&\,\,\,\,\,\,\,+\,\,c_{3}(n,p)\left|h\right|^{2}\Vert Df_{\varepsilon}\Vert_{L^{p'}(Q_{R_{0}})}\,\Vert Du_{\varepsilon}\Vert_{L^{p}(Q_{R_{0}})},
\end{align*}which holds for every $t_{1}\in(t_{0}-\rho^{2},t_{0})$ and every
suitably small $h\neq0$. \\
We now choose a cut-off function $\eta\in C_{0}^{\infty}\left(B_{\rho}(x_{0})\right)$
with $\eta\equiv1$ on $B_{\rho/2}(x_{0})$ such that $0\leq\eta\leq1$
and $\left|D\eta\right|\leq C/\rho$. For the cut-off function in
time, we choose the piecewise affine function $\chi:(T_{1},T_{2})\rightarrow\left[0,1\right]$
with $\chi\equiv0$ on $(T_{1},t_{0}-\rho^{2})$, $\chi\equiv1$ on
$(t_{0}-(\rho/2)^{2},T_{2})$ and $\partial_{t}\chi\equiv\frac{4}{3\rho^{2}}$
on $(t_{0}-\rho^{2},t_{0}-(\rho/2)^{2})$.\\
Dividing both sides of the previous estimate by $\left|h\right|^{2}$
and using the properties of $\chi$ and $\eta$, we obtain\begin{equation}\label{eq:rif8}
\begin{split}
&\sup_{t_{0}-(\rho/2)^{2}<t<t_{0}}\int_{B_{\rho/2}(x_{0})}\left|\Delta_{h}u_{\varepsilon}(x,t)\right|^{2}dx\,+\int_{Q_{\rho/2}(z_{0})}\left|\Delta_{h}H_{\frac{p}{2}}(Du_{\varepsilon})\right|^{2}dz\\
&\,\,\,\,\,\,\,\,\,\,\leq\,c_{4}(n,p)\,\rho^{-2}\int_{Q_{\rho}(z_{0})}\left[\left(\left|Du_{\varepsilon}(x+he_{s},t)\right|^{\frac{p-2}{2}}+\left|Du_{\varepsilon}\right|^{\frac{p-2}{2}}\right)^{2}+1\right]\left|\Delta_{h}u_{\varepsilon}\right|^{2}dz\\
&\,\,\,\,\,\,\,\,\,\,\,\,\,\,\,\,\,+\,c_{4}(n,p)\,\Vert Df_{\varepsilon}\Vert_{L^{p'}(Q_{R_{0}})}\,\Vert Du_{\varepsilon}\Vert_{L^{p}(Q_{R_{0}})}.
\end{split}
\end{equation}Now we set
\[
I_{5}:=\int_{Q_{\rho}(z_{0})}\left[\left(\left|Du_{\varepsilon}(x+he_{s},t)\right|^{\frac{p-2}{2}}+\left|Du_{\varepsilon}\right|^{\frac{p-2}{2}}\right)^{2}+1\right]\left|\Delta_{h}u_{\varepsilon}\right|^{2}dz.
\]
$\hspace*{1em}$If $p>2$, using Hölder's inequality with exponents
$\left(\frac{p}{2},\frac{p}{p-2}\right)$ and the properties of the
difference quotients, we can control $I_{5}$ as follows\begin{equation}\label{eq:rif9}
\begin{split}
I_{5}\,&\leq\, c_{5}(n)\left(\int_{Q_{R_{0}}}\left|Du_{\varepsilon}\right|^{p}dz\right)^{\frac{2}{p}}\left(\int_{Q_{R_{0}}}\left[\left|Du_{\varepsilon}\right|^{p-2}+1\right]^{\frac{p}{p-2}}dz\right)^{\frac{p-2}{p}}\\
&\leq\, c_{6}(n,p)\left(\int_{Q_{R_{0}}}\left|Du_{\varepsilon}\right|^{p}dz\right)^{\frac{2}{p}}\left(\int_{Q_{R_{0}}}\left[\left|Du_{\varepsilon}\right|^{p}+1\right]\,dz\right)^{\frac{p-2}{p}}\\
&\leq\,c_{6}(n,p)\,\int_{Q_{R_{0}}}\left(\left|Du_{\varepsilon}\right|^{p}+1\right)\,dz,
\end{split}
\end{equation}provided that $\left|h\right|$ is sufficiently small. Joining estimates
\eqref{eq:rif8} and \eqref{eq:rif9}, for $p>2$ we then have\begin{align*}
&\sup_{t_{0}-(\rho/2)^{2}<t<t_{0}}\int_{B_{\rho/2}(x_{0})}\left|\Delta_{s,h}u_{\varepsilon}(x,t)\right|^{2}dx\,+\int_{Q_{\rho/2}(z_{0})}\left|\Delta_{s,h}H_{\frac{p}{2}}(Du_{\varepsilon})\right|^{2}dz\\
&\,\,\,\,\,\,\,\,\,\,\,\,\,\,\leq c\left[\Vert Df_{\varepsilon}\Vert_{L^{p'}(Q_{R_{0}})}\left(\int_{Q_{R_{0}}}\left|Du_{\varepsilon}\right|^{p}dz\right)^{\frac{1}{p}}+\,\rho^{-2}\int_{Q_{R_{0}}}\left(\left|Du_{\varepsilon}\right|^{p}+1\right)\,dz\right],
\end{align*}with $c\equiv c(n,p)>0$. Since the above estimate holds for every
$s\in\left\{ 1,\ldots,n\right\} $ and every sufficiently small $h\in\mathbb{R}\setminus\left\{ 0\right\} $,
for $p>2$, by Lemma \ref{lem:RappIncre} we may conclude that
\[
H_{\frac{p}{2}}(Du_{\varepsilon})\,\in\,L_{loc}^{2}\left(T_{1},T_{2};W_{loc}^{1,2}(\Omega',\mathbb{R}^{n})\right).
\]
 $\hspace*{1em}$Finally, if $p=2$, arguing in a similar fashion
we reach the same conclusions. Moreover, letting $h\rightarrow0$
in the above estimate, we obtain the Caccioppoli-type inequality \eqref{eq:Caccio1},
which in turn implies the validity of (\ref{eq:apriori1}).\end{proof}\newpage{}

\noindent $\hspace*{1em}$As a consequence of Proposition \ref{prop:apriori1},
we are able to establish the following higher integrability result
for the spatial gradient $Du_{\varepsilon}$, whose proof is based
on Sobolev embedding theorem in space applied slicewise: 
\end{singlespace}
\begin{prop}
\begin{singlespace}
\noindent \label{prop:apriori2} Under the assumptions of Proposition
\ref{prop:apriori1}, we obtain that 
\[
Du_{\varepsilon}\,\in\,L_{loc}^{p\,+\,\frac{4}{n}}\left(\Omega'_{1,2},\mathbb{R}^{n}\right).
\]
Moreover, there exists a positive constant $C_{1}=C_{1}(n,p)$ such
that for any cylinder $Q_{\gamma}(z_{0})\subset Q_{\rho}(z_{0})\Subset\Omega'_{1,2}$
there holds\begin{align*}
\int_{Q_{\gamma}(z_{0})}(\vert Du_{\varepsilon}\vert-\nu)_{+}^{p\,+\,\frac{4}{n}}\,\,dz\,\,\leq\,\,C_{1}&\left[\sup_{t_{0}-\rho^{2}<t<t_{0}}\int_{B_{\rho}(x_{0})}\left|Du_{\varepsilon}(x,t)\right|^{2}dx\right]^{\frac{2}{n}}\cdot\\
&\cdot\int_{Q_{\rho}(z_{0})}\left(\left|DH_{\frac{p}{2}}(Du_{\varepsilon})\right|^{2}+\frac{1}{(\rho-\gamma)^{2}}\left|Du_{\varepsilon}\right|^{p}\right)dz.
\end{align*}In addition, the following estimate \begin{equation}\label{eq:apriori2}
\begin{split}
&\int_{Q_{\rho/2}(z_{0})}(\vert Du_{\varepsilon}\vert-\nu)_{+}^{p\,+\,\frac{4}{n}}\,\,dz\\
&\,\,\,\,\,\,\,\,\,\,\,\,\,\,\leq\,C_{2}\left[\Vert Df_{\varepsilon}\Vert_{L^{p'}(Q_{R_{0}})}\left(\int_{Q_{R_{0}}}\left|Du_{\varepsilon}\right|^{p}dz\right)^{\frac{1}{p}}+\,\rho^{-2}\int_{Q_{R_{0}}}\left(\left|Du_{\varepsilon}\right|^{p}+1\right)\,dz\right]^{1\,+\,\frac{2}{n}}
\end{split}
\end{equation}holds true for any parabolic cylinder $Q_{\rho}(z_{0})\subset Q_{2\rho}(z_{0})\subset Q_{R_{0}}(z_{0})\Subset\Omega'_{1,2}$
and a positive constant $C_{2}=C_{2}(n,p)$.
\end{singlespace}
\end{prop}

\begin{singlespace}
\noindent \begin{proof}[\bfseries{Proof}] Let us introduce the following
function\medskip{}
\[
G_{\varepsilon}:=\,\frac{2\,n}{np+4}\,\left|H_{\frac{p}{2}}(Du_{\varepsilon})\right|^{\frac{4}{np}\,+\,1}.
\]
\\
In what follows, we will denote by $c_{k}$ some positive constants
which do not depend on $\varepsilon$. From the definition of $G_{\varepsilon}$,
we obtain\\
\begin{align}
\left|DG_{\varepsilon}\right| & \,=\,\frac{2}{p}\,\left|H_{\frac{p}{2}}(Du_{\varepsilon})\right|^{\frac{4}{np}}\left|D\left|H_{\frac{p}{2}}(Du_{\varepsilon})\right|\right|\nonumber \\
 & \,\leq\,c_{1}\,\left|H_{\frac{p}{2}}(Du_{\varepsilon})\right|^{\frac{4}{np}}\left|D\left[H_{\frac{p}{2}}(Du_{\varepsilon})\right]\right|,\label{eq:DV-1}
\end{align}
\\
where $c_{1}\equiv c_{1}(n,p)>0$. For $B_{\rho}(x_{0})\Subset\Omega'$,
let $\varphi\in C_{0}^{\infty}\left(B_{\rho}(x_{0})\right)$ and $\chi\in W^{1,\infty}\left((T_{1},T_{2})\right)$
be two non-negative cut-off functions with $\chi(T_{1})=0$ and $\partial_{t}\chi\geq0$.
Now, fix a time $t_{0}\in(T_{1},T_{2})$ and apply the Sobolev embedding
theorem on the time slices $\varSigma_{t}:=B_{\rho}(x_{0})\times\left\{ t\right\} $
for almost every $t\in(T_{1},t_{0})$, to infer that\\
\begin{align*}
\int_{\varSigma_{t}}\varphi^{2}G_{\varepsilon}^{2}\,dx\,\,&\leq \,c_{2}(n)\left(\int_{\varSigma_{t}}\left|D(\varphi \,G_{\varepsilon})\right|^{\frac{2n}{n+2}}\,dx\right)^{\frac{n+2}{n}}\\
&=\,c_{2}(n)\left(\int_{\varSigma_{t}}\left|\varphi\,DG_{\varepsilon}+G_{\varepsilon}\,D\varphi\right|^{\frac{2n}{n+2}}\,dx\right)^{\frac{n+2}{n}}\\
&\leq \,c_{3}(n)\left(\int_{\varSigma_{t}}\left|\varphi\,DG_{\varepsilon}\right|^{\frac{2n}{n+2}}\,dx\right)^{\frac{n+2}{n}}+\,c_{3}(n)\left(\int_{\varSigma_{t}}\left|G_{\varepsilon}\,D\varphi\right|^{\frac{2n}{n+2}}\,dx\right)^{\frac{n+2}{n}}\\
&=:\,c_{3}\,J_{1}(t)\,+\,c_{3}\,J_{2}(t),
\end{align*}where, in the second to last line, we have applied Minkowski's and
Young's inequalities one after the other. In the following, we estimate
$J_{1}(t)$ and $J_{2}(t)$ separately. Let us first consider $J_{1}(t)$.
Using (\ref{eq:DV-1}) and Hölder's inequality, we deduce\begin{align*}
J_{1}(t)\,\,&\leq \,c_{4}(n,p)\left(\int_{\varSigma_{t}}\varphi^{\frac{2n}{n+2}}\left[\left(\left|Du_{\varepsilon}\right|-\nu\right)_{+}^{2/n}\left|DH_{\frac{p}{2}}(Du_{\varepsilon})\right|\right]^{\frac{2n}{n+2}}dx\right)^{\frac{n+2}{n}}\\
&\leq \,c_{4}(n,p)\int_{\varSigma_{t}}\varphi^{2}\left|DH_{\frac{p}{2}}(Du_{\varepsilon})\right|^{2}dx\,\left(\int_{\mathrm{supp}(\varphi)}\left(\left|Du_{\varepsilon}\right|-\nu\right)_{+}^{2}dx\right)^{\frac{2}{n}}\\
&\leq \,c_{4}(n,p)\int_{\varSigma_{t}}\varphi^{2}\left|DH_{\frac{p}{2}}(Du_{\varepsilon})\right|^{2}dx\,\left(\int_{\mathrm{supp}(\varphi)}\left|Du_{\varepsilon}\right|^{2}dx\right)^{\frac{2}{n}}.
\end{align*}\\
We now turn our attention to $J_{2}(t)$. Using the definition of
$G_{\varepsilon}$ and Hölder's inequality, we can conclude\begin{align*}
J_{2}(t)\,&=\,\frac{4\,n^{2}}{(np+4)^{2}}\left(\int_{\varSigma_{t}}\left(\left|Du_{\varepsilon}\right|-\nu\right)_{+}^{\frac{np\,+\,4}{n+2}}\left|D\varphi\right|^{\frac{2n}{n+2}}\,dx\right)^{\frac{n+2}{n}}\\
&\leq\,\frac{4\,n^{2}}{(np+4)^{2}}\left(\int_{\varSigma_{t}}\left[\left|D\varphi\right|^{2}\left|Du_{\varepsilon}\right|^{p}\right]^{\frac{n}{n+2}}\left|Du_{\varepsilon}\right|^{\frac{4}{n+2}}\,dx\right)^{\frac{n+2}{n}}\\
&\leq\,\frac{4\,n^{2}}{(np+4)^{2}}\int_{\varSigma_{t}}\left|D\varphi\right|^{2}\left|Du_{\varepsilon}\right|^{p}dx\,\left(\int_{\mathrm{supp}(\varphi)}\left|Du_{\varepsilon}\right|^{2}dx\right)^{\frac{2}{n}}.
\end{align*}\\
Putting together the last three estimates and integrating with respect
to time, we obtain\\
\begin{equation}\label{eq:rif10}
\begin{split}
\int_{Q_{T_{1},t_{0}}}\chi\varphi^{2}\,(\vert Du_{\varepsilon}\vert-\nu)_{+}^{p\,+\,\frac{4}{n}}\,\,dz\,&\leq \,\,c_{5}(n,p)\int_{T_{1}}^{t_{0}}\chi\left[\int_{\mathrm{supp}(\varphi)}\left|Du_{\varepsilon}(x,t)\right|^{2}dx\right]^{\frac{2}{n}}\cdot\\
&\,\,\,\,\,\,\,\,\,\,\,\,\,\,\,\,\,\,\,\,\,\,\,\,\,\,\cdot\left[\int_{\varSigma_{t}}\left(\varphi^{2}\left|DH_{\frac{p}{2}}(Du_{\varepsilon})\right|^{2}+\left|D\varphi\right|^{2}\left|Du_{\varepsilon}\right|^{p}\right)dx\right]dt\\
&\leq \,\,c_{5}(n,p)\left[\sup_{T_{1}<t<t_{0},\,\chi(t)\neq0}\int_{\mathrm{supp}(\varphi)}\left|Du_{\varepsilon}(x,t)\right|^{2}dx\right]^{\frac{2}{n}}\cdot\\
&\,\,\,\,\,\,\,\,\,\,\,\,\,\,\,\,\,\,\,\,\,\,\,\,\,\,\cdot\int_{Q_{T_{1},t_{0}}}\chi\left(\varphi^{2}\left|DH_{\frac{p}{2}}(Du_{\varepsilon})\right|^{2}+\left|D\varphi\right|^{2}\left|Du_{\varepsilon}\right|^{p}\right)dz,
\end{split}
\end{equation}where we have used the abbreviation $Q_{T_{1},t_{0}}:=B_{\rho}(x_{0})\times(T_{1},t_{0})$.
Now we perform a particular choice of the cut-off functions $\chi$
and $\varphi$ involved above. For a parabolic cylinder $Q_{\rho}(z_{0})\Subset\Omega'_{1,2}$
we choose $\chi\in W^{1,\infty}\left((T_{1},T_{2})\right)$ such that
\[
\chi\equiv0\,\,\,\,\,\mathrm{on}\,\,\left(T_{1},t_{0}-\rho^{2}\right],\,\,\,\,\,\,\,\,\,\chi\equiv1\,\,\,\,\,\mathrm{on}\,\,\left[t_{0}-\gamma^{2},T_{2}\right)\,\,\,\,\,\,\,\,\,\mathrm{and}\,\,\,\,\,\,\,\,\,\partial_{t}\chi\geq0,
\]
with $0<\gamma<\rho$. As for $\varphi\in C_{0}^{\infty}\left(B_{\rho}(x_{0})\right)$,
we assume that $\varphi\equiv1$ on $B_{\gamma}(x_{0})$, $0\leq\varphi\leq1$
and $\left|D\varphi\right|\leq\frac{C}{\rho-\gamma}$. With these
choices \eqref{eq:rif10} turns into\begin{equation}\label{eq:rif11}
\begin{split}
\int_{Q_{\gamma}(z_{0})}(\vert Du_{\varepsilon}\vert-\nu)_{+}^{p\,+\,\frac{4}{n}}\,\,dz\,\,\leq\,\,c_{6}(n,p)&\left[\sup_{t_{0}-\rho^{2}<t<t_{0}}\int_{B_{\rho}(x_{0})}\left|Du_{\varepsilon}(x,t)\right|^{2}dx\right]^{\frac{2}{n}}\cdot\\
&\cdot\int_{Q_{\rho}(z_{0})}\left(\left|DH_{\frac{p}{2}}(Du_{\varepsilon})\right|^{2}+\frac{1}{(\rho-\gamma)^{2}}\left|Du_{\varepsilon}\right|^{p}\right)dz.
\end{split}
\end{equation}We now choose $\gamma=\rho/2$ and use \eqref{eq:Caccio1} with $Q_{\rho}(z_{0})$
replaced by $Q_{2\rho}(z_{0})$, in order to estimate the first and
second integral on the right-hand side of \eqref{eq:rif11}. After
changing notation about the cylinders involved, we finally obtain
the inequality \eqref{eq:apriori2}, which ensures that $Du_{\varepsilon}\in L_{loc}^{p\,+\,\frac{4}{n}}\left(\Omega'_{1,2},\mathbb{R}^{n}\right)$.\end{proof}
\end{singlespace}
\begin{singlespace}

\section{Proof of Theorem \ref{thm:main4} \label{sec:main theo} }
\end{singlespace}

\begin{singlespace}
\noindent $\hspace*{1em}$We now prove Theorem \ref{thm:main4}, by
dividing the proof into three steps. The first step consists in constructing
a family of Cauchy-Dirichlet problems, for which we are allowed to
use the \textit{a priori} estimates from Propositions \ref{prop:apriori1}
and \ref{prop:apriori2}. At this stage, Lemma \ref{lem:giusti} will
play a key role in deriving an \textit{a priori} estimate for the
solutions to these problems (i.e. the comparison maps): we specifically
refer to estimate \eqref{eq:essenz1} below. \\
$\hspace*{1em}$In the second step, we will show that the $L^{p}$-norms
of the spatial gradients of the comparison maps are actually uniformly
bounded (see estimate \eqref{eq:ee18}), and this is where Lemma \ref{lem:inter}
will come into play.\\
$\hspace*{1em}$Finally, in the third step, we shall use a standard
comparison argument, as well as the results obtained in the previous
steps, to reach the desired conclusion.\\

\noindent \begin{proof}[\bfseries{Proof of Theorem~\ref{thm:main4}}]
\textbf{Step 1: }\textbf{\textit{a priori}}\textbf{ estimate for the
comparison maps.}\\
\\
We shall keep the notation introduced for the proof of Proposition
\ref{prop:apriori1}, starting from the case $p>2$. For $\varepsilon\in(0,1]$
and a couple of standard, non-negative, radially symmetric mollifiers
$\phi_{1}\in C_{0}^{\infty}(B_{1}(0))$ and $\phi_{2}\in C_{0}^{\infty}((-1,1))$,
we define the function $f_{\varepsilon}$ as in (\ref{eq:molli}),
where $f$ is meant to be extended by zero outside $\Omega_{T}$.
Therefore, we have that $f_{\varepsilon}\in C^{\infty}(\Omega_{T})$.
\\
$\hspace*{1em}$Now, for any fixed $\varepsilon\in(0,1]$ let us define
the comparison map\vspace{0.4cm}
\[
u_{\varepsilon}\in C^{0}\left([t_{0}-R_{0}^{2},t_{0}];L^{2}\left(B_{R_{0}}(x_{0})\right)\right)\cap L^{p}\left(t_{0}-R_{0}^{2},t_{0};W^{1,p}\left(B_{R_{0}}(x_{0})\right)\right)
\]
as the unique energy solution of the Cauchy-Dirichlet problem\vspace{0.4cm}
\begin{equation}
\begin{cases}
\begin{array}{cc}
\partial_{t}u_{\varepsilon}-\mathrm{div}\left(H_{p-1}(Du_{\varepsilon})+\varepsilon\,\vert Du_{\varepsilon}\vert^{p-2}Du_{\varepsilon}\right)=f_{\varepsilon} & \,\,\,\mathrm{in}\,\,\,Q_{R_{0}}(z_{0})\\
u_{\varepsilon}=u & \,\,\,\,\,\,\,\,\,\,\,\,\,\,\mathrm{on}\,\,\,\partial_{\mathrm{par}}Q_{R_{0}}(z_{0}),
\end{array}\end{cases}\label{eq:Diri}
\end{equation}
\\
where $Q_{R_{0}}(z_{0})\Subset\Omega_{T}$ and the initial-boundary
condition is meant in the sense of Definition \ref{def:L2-traces}
(see \cite[Chapter 2]{Lions} or \cite[Chapter 9]{DiBene} for the
existence). Moreover, let us fix a positive number $R<R_{0}$ and
arbitrary radii
\[
\frac{R}{2}\leq r<\ell<\rho<\gamma<\lambda r<R,
\]
\\
with $1<\lambda<2$. In what follows, we will denote by $c_{k}$ some
positive constants which do not depend either on $h$ or $\varepsilon$.\\
For a fixed time $t_{1}\in(t_{0}-\rho^{2},t_{0})$ and $\delta\in(0,t_{0}-t_{1})$,
we choose $\Phi(x,t)=\widetilde{\chi}(t)\chi(t)\eta^{2}(x)$ with
$\chi\in W^{1,\infty}\left((t_{0}-R_{0}^{2},t_{0}),\left[0,1\right]\right)$,
$\chi\equiv0$ on $(t_{0}-R_{0}^{2},t_{0}-\rho^{2})$ and $\partial_{t}\chi\geq0$,
$\eta\in C_{0}^{\infty}\left(B_{\rho}(x_{0}),\left[0,1\right]\right)$,
and with the Lipschitz continuous function $\widetilde{\chi}:(t_{0}-R_{0}^{2},t_{0})\rightarrow\mathbb{R}$
defined by\\
\\
\[
\widetilde{\chi}(t)=\begin{cases}
\begin{array}{c}
1\\
\mathrm{affine}\\
0
\end{array} & \begin{array}{c}
\mathrm{if}\,\,t\leq t_{1},\quad\quad\quad\;\,\,\\
\mathrm{if}\,\,t_{1}<t<t_{1}+\delta,\\
\mathrm{if}\,\,t\geq t_{1}+\delta.\quad\quad
\end{array}\end{cases}
\]
\\
\\
Setting $Q^{t_{1}}:=B_{\rho}(x_{0})\times(t_{0}-\rho^{2},t_{1})$
and arguing as in the first part of the proof of Proposition \ref{prop:apriori1},
from $(\ref{eq:Diri})_{1}$ we get\vspace{0.4cm}
\begin{equation}\label{eq:ee1}
\begin{split}
&\frac{1}{2}\int_{B_{\rho}(x_{0})}\chi(t_{1})\eta^{2}(x)\left|\tau_{h}u_{\varepsilon}(x,t_{1})\right|^{2}dx\,+\int_{Q^{t_{1}}}\chi(t)\eta^{2}(x)\langle\tau_{h}H_{p-1}(Du_{\varepsilon}),D\tau_{h}u_{\varepsilon}\rangle\,dz\\
&+\,\varepsilon\int_{Q^{t_{1}}}\chi(t)\eta^{2}(x)\langle\tau_{h}\left[\vert Du_{\varepsilon}\vert^{p-2}Du_{\varepsilon}\right],D\tau_{h}u_{\varepsilon}\rangle\,dz\\
&=-2\int_{Q^{t_{1}}}\chi(t)\eta(x)\langle\tau_{h}H_{p-1}(Du_{\varepsilon}),D\eta\rangle\tau_{h}u_{\varepsilon}\,\,dz\,+\int_{Q^{t_{1}}}(\tau_{h}f_{\varepsilon})(\tau_{h}u_{\varepsilon})\chi(t)\eta^{2}(x)\,dz\\
&\,\,\,\,\,\,-2\,\varepsilon\int_{Q^{t_{1}}}\chi(t)\eta(x)\langle\tau_{h}\left[\vert Du_{\varepsilon}\vert^{p-2}Du_{\varepsilon}\right],D\eta\rangle\tau_{h}u_{\varepsilon}\,\,dz\,+\,\frac{1}{2}\int_{Q^{t_{1}}}(\partial_{t}\chi)\eta^{2}(x)\left|\tau_{h}u_{\varepsilon}\right|^{2}dz\\
&=:A_{1}+A_{2}+A_{3}+A_{4},
\end{split}
\end{equation}for every sufficiently small $h\in\mathbb{R}\setminus\left\{ 0\right\} $
and every $t_{1}\in(t_{0}-\rho^{2},t_{0})$. In the following, we
estimate $A_{1}$, $A_{2}$ and $A_{3}$ separately. Let us first
consider $A_{1}$. Using Lemma \ref{lem:Brasco} together with Young's
inequality with $\sigma>0$ and exponents $\left(2,2\right)$, we
obtain\vspace{0.4cm}
\begin{equation}\label{eq:ee2}
\begin{split}
\left|A_{1}\right|\leq\,&\,\,\frac{2(p-1)^{2}}{\sigma}\int_{Q^{t_{1}}}\chi(t)\left|D\eta\right|^{2}\left(\left|H_{\frac{p}{2}}(Du_{\varepsilon}(x+he_{s},t))\right|^{\frac{p-2}{p}}+\left|H_{\frac{p}{2}}(Du_{\varepsilon})\right|^{\frac{p-2}{p}}\right)^{2}\left|\tau_{h}u_{\varepsilon}\right|^{2}dz\\
&+\,\frac{\sigma}{2}\int_{Q^{t_{1}}}\chi(t)\eta^{2}(x)\left|\tau_{h}H_{\frac{p}{2}}(Du_{\varepsilon})\right|^{2}dz.
\end{split}
\end{equation}\newpage As for $A_{3}$, by Lemma \ref{lem:Lind} we similarly have\vspace{0.4cm}
\begin{equation}\label{eq:ee3}
\begin{split}
\left|A_{3}\right|\leq\,&\,\,\frac{2\varepsilon\,(p-1)^{2}}{\sigma}\int_{Q^{t_{1}}}\chi(t)\left|D\eta\right|^{2}\left(\left|Du_{\varepsilon}(x+he_{s},t)\right|^{\frac{p-2}{2}}+\left|Du_{\varepsilon}\right|^{\frac{p-2}{2}}\right)^{2}\left|\tau_{h}u_{\varepsilon}\right|^{2}dz\\
&+\,\frac{\varepsilon\sigma}{2}\int_{Q^{t_{1}}}\chi(t)\eta^{2}(x)\left|\tau_{h}V_{p}(Du_{\varepsilon})\right|^{2}dz.
\end{split}
\end{equation}\\
We now turn our attention to $A_{2}$. Thanks to Propositions \ref{prop:apriori1}
and \ref{prop:apriori2}, we know that\vspace{0.4cm}
\[
H_{\frac{p}{2}}(Du_{\varepsilon})\,\in\,L_{loc}^{2}\left(t_{0}-R_{0}^{2},t_{0};W_{loc}^{1,2}\left(B_{R_{0}}(x_{0}),\mathbb{R}^{n}\right)\right),
\]
\vspace{0.1cm}
\[
Du_{\varepsilon}\,\in\,L_{loc}^{\infty}\left(t_{0}-R_{0}^{2},t_{0};L_{loc}^{2}\left(B_{R_{0}}(x_{0}),\mathbb{R}^{n}\right)\right)\cap L_{loc}^{p\,+\,\frac{4}{n}}\left(Q_{R_{0}}(z_{0}),\mathbb{R}^{n}\right)
\]
\\
and\\
\begin{equation}\label{eq:ee3bis}
\begin{split}
\int_{Q_{\gamma}(z_{0})}(\vert Du_{\varepsilon}\vert-\nu)_{+}^{p\,+\,\frac{4}{n}}\,\,dz\,\,\leq\,\,c_{1}(n,p)&\left[\sup_{t_{0}-(\lambda r)^{2}<t<t_{0}}\int_{B_{\lambda r}(x_{0})}\left|Du_{\varepsilon}(x,t)\right|^{2}dx\right]^{\frac{2}{n}}\cdot\\
&\cdot\int_{Q_{\lambda r}(z_{0})}\left(\left|DH_{\frac{p}{2}}(Du_{\varepsilon})\right|^{2}+\frac{1}{(\lambda r-\gamma)^{2}}\left|Du_{\varepsilon}\right|^{p}\right)dz.
\end{split}
\end{equation}\\
\\
Therefore, taking advantage of the properties of $f_{\varepsilon}$,
$u_{\varepsilon}$, $\chi$ and $\eta$, and using Hölder's inequality
with exponents $\left(\frac{np\,+\,4}{np\,+\,4\,-\,n},p+\frac{4}{n}\right)$
together with Lemma \ref{lem:Giusti1}, we obtain\vspace{0.4cm}
\begin{equation}
A_{2}\leq\,c_{2}(n,p)\left|h\right|^{2}\left(\int_{Q_{\gamma}(z_{0})}\left|Df_{\varepsilon}\right|^{\frac{np\,+\,4}{np\,+\,4\,-\,n}}dz\right)^{\frac{np\,+\,4\,-\,n}{np\,+\,4}}\left(\int_{Q_{\gamma}\left(z_{0}\right)}\left|Du_{\varepsilon}\right|^{p\,+\,\frac{4}{n}}dz\right)^{\frac{n}{np\,+\,4}},\label{eq:ee4}
\end{equation}
\\
provided that $\left|h\right|$ is suitably small. Now, by virtue
of Lemma \ref{lem:Brasco} we have\vspace{0.4cm}
\begin{equation}
\frac{4}{p^{2}}\int_{Q^{t_{1}}}\chi(t)\eta^{2}(x)\left|\tau_{h}H_{\frac{p}{2}}(Du_{\varepsilon})\right|^{2}dz\,\leq\,\int_{Q^{t_{1}}}\chi(t)\eta^{2}(x)\langle\tau_{h}H_{p-1}(Du_{\varepsilon}),D\tau_{h}u_{\varepsilon}\rangle\,dz.\label{eq:ee5}
\end{equation}
\\
Similarly, by Lemma \ref{lem:Lind} we obtain\vspace{0.4cm}
\begin{equation}
\frac{4\,\varepsilon}{p^{2}}\int_{Q^{t_{1}}}\chi(t)\eta^{2}(x)\left|\tau_{h}V_{p}(Du_{\varepsilon})\right|^{2}dz\,\leq\,\varepsilon\int_{Q^{t_{1}}}\chi(t)\eta^{2}(x)\langle\tau_{h}\left[\vert Du_{\varepsilon}\vert^{p-2}Du_{\varepsilon}\right],D\tau_{h}u_{\varepsilon}\rangle\,dz.\label{eq:ee6}
\end{equation}
\\
Collecting estimates \eqref{eq:ee1}, \eqref{eq:ee2}, \eqref{eq:ee3},
(\ref{eq:ee4}), (\ref{eq:ee5}) and (\ref{eq:ee6}), we arrive at\vspace{0.4cm}
\begin{equation}\label{eq:ee7}
\begin{split}
&\int_{B_{\rho}(x_{0})}\chi(t_{1})\eta^{2}(x)\left|\tau_{h}u_{\varepsilon}(x,t_{1})\right|^{2}dx\,+\int_{Q^{t_{1}}}\chi(t)\eta^{2}(x)\left|\tau_{h}H_{\frac{p}{2}}(Du_{\varepsilon})\right|^{2}dz\\
&+\,\,\varepsilon\int_{Q^{t_{1}}}\chi(t)\eta^{2}(x)\left|\tau_{h}V_{p}(Du_{\varepsilon})\right|^{2}dz\\
&\leq\,\sigma\,c_{3}(p)\int_{Q^{t_{1}}}\chi(t)\eta^{2}(x)\left|\tau_{h}H_{\frac{p}{2}}(Du_{\varepsilon})\right|^{2}dz\,+\,\varepsilon\,\sigma\,c_{3}(p)\int_{Q^{t_{1}}}\chi(t)\eta^{2}(x)\left|\tau_{h}V_{p}(Du_{\varepsilon})\right|^{2}dz\\
&\,\,\,\,\,\,\,+\,\,\frac{c_{3}(p)}{\sigma}\int_{Q^{t_{1}}}\chi(t)\left|D\eta\right|^{2}\left(\left|H_{\frac{p}{2}}(Du_{\varepsilon}(x+he_{s},t))\right|^{\frac{p-2}{p}}+\left|H_{\frac{p}{2}}(Du_{\varepsilon})\right|^{\frac{p-2}{p}}\right)^{2}\left|\tau_{h}u_{\varepsilon}\right|^{2}dz\\
&\,\,\,\,\,\,\,+\,\,\frac{\varepsilon\,c_{3}(p)}{\sigma}\int_{Q^{t_{1}}}\chi(t)\left|D\eta\right|^{2}\left(\left|Du_{\varepsilon}(x+he_{s},t)\right|^{\frac{p-2}{2}}+\left|Du_{\varepsilon}\right|^{\frac{p-2}{2}}\right)^{2}\left|\tau_{h}u_{\varepsilon}\right|^{2}dz\\
&\,\,\,\,\,\,\,+\,\,c_{4}(n,p)\left|h\right|^{2}\left(\int_{Q_{\gamma}(z_{0})}\left|Df_{\varepsilon}\right|^{\frac{np\,+\,4}{np\,+\,4\,-\,n}}dz\right)^{\frac{np\,+\,4\,-\,n}{np\,+\,4}}\left(\int_{Q_{\gamma}\left(z_{0}\right)}\left|Du_{\varepsilon}\right|^{p\,+\,\frac{4}{n}}dz\right)^{\frac{n}{np\,+\,4}}\\
&\,\,\,\,\,\,\,+\,\,c_{3}(p)\int_{Q^{t_{1}}}(\partial_{t}\chi)\eta^{2}(x)\left|\tau_{h}u_{\varepsilon}\right|^{2}dz.
\end{split}
\end{equation}\\
Choosing $\sigma=(2\,c_{3}(p))^{-1}$ and reabsorbing the first two
integrals in the right-hand side of \eqref{eq:ee7} by the left-hand
side, we get\begin{align*}
&\int_{B_{\rho}(x_{0})}\chi(t_{1})\eta^{2}(x)\left|\tau_{h}u_{\varepsilon}(x,t_{1})\right|^{2}dx\,+\int_{Q^{t_{1}}}\chi(t)\eta^{2}(x)\left|\tau_{h}H_{\frac{p}{2}}(Du_{\varepsilon})\right|^{2}dz\\
&+\,\,\varepsilon\int_{Q^{t_{1}}}\chi(t)\eta^{2}(x)\left|\tau_{h}V_{p}(Du_{\varepsilon})\right|^{2}dz\\
&\leq\,c_{5}(p)\int_{Q^{t_{1}}}\chi(t)\left|D\eta\right|^{2}\left(\left|H_{\frac{p}{2}}(Du_{\varepsilon}(x+he_{s},t))\right|^{\frac{p-2}{p}}+\left|H_{\frac{p}{2}}(Du_{\varepsilon})\right|^{\frac{p-2}{p}}\right)^{2}\left|\tau_{h}u_{\varepsilon}\right|^{2}dz\\
&\,\,\,\,\,\,\,+\,\,\varepsilon\,c_{5}(p)\int_{Q^{t_{1}}}\chi(t)\left|D\eta\right|^{2}\left(\left|Du_{\varepsilon}(x+he_{s},t)\right|^{\frac{p-2}{2}}+\left|Du_{\varepsilon}\right|^{\frac{p-2}{2}}\right)^{2}\left|\tau_{h}u_{\varepsilon}\right|^{2}dz\\
&\,\,\,\,\,\,\,+\,\,c_{6}(n,p)\left|h\right|^{2}\left(\int_{Q_{\gamma}(z_{0})}\left|Df_{\varepsilon}\right|^{\frac{np\,+\,4}{np\,+\,4\,-\,n}}dz\right)^{\frac{np\,+\,4\,-\,n}{np\,+\,4}}\left(\int_{Q_{\gamma}\left(z_{0}\right)}\left|Du_{\varepsilon}\right|^{p\,+\,\frac{4}{n}}dz\right)^{\frac{n}{np\,+\,4}}\\
&\,\,\,\,\,\,\,+\,\,c_{5}(p)\int_{Q^{t_{1}}}(\partial_{t}\chi)\eta^{2}(x)\left|\tau_{h}u_{\varepsilon}\right|^{2}dz,
\end{align*}which holds for every $t_{1}\in(t_{0}-\rho^{2},t_{0})$ and every
sufficiently small $h\in\mathbb{R}\setminus\left\{ 0\right\} $.\\
We now choose a cut-off function $\eta\in C_{0}^{\infty}\left(B_{\rho}(x_{0})\right)$
with $\eta\equiv1$ on $B_{\ell}(x_{0})$ such that $0\leq\eta\leq1$
and $\left|D\eta\right|\leq C/(\rho-\ell)$. For the cut-off function
in time, we choose $\chi\in W^{1,\infty}\left((t_{0}-R_{0}^{2},t_{0}),\left[0,1\right]\right)$
such that
\[
\chi\equiv0\,\,\,\,\,\,\,\,\mathrm{on}\,\,\left(t_{0}-R_{0}^{2},t_{0}-\rho^{2}\right],
\]
\vspace{0.2cm}
\[
\chi\equiv1\,\,\,\,\,\,\,\,\mathrm{on}\,\,\left[t_{0}-\ell^{2},t_{0}\right)\,\,\,\,\,\,\,\,\,\,\,\,\,\,\,\,\,
\]
\vspace{0.2cm}
and\vspace{0.2cm}
\[
\partial_{t}\chi\leq\frac{C}{(\rho-\ell)^{2}}\,\,\,\,\,\,\,\,\mathrm{on}\,\,\left(t_{0}-\rho^{2},t_{0}-\ell^{2}\right).\,\,\,\,\,\,\,\,\,\,\,\,\,\,\,\,\,\,\,\,\,\,\,\,\,
\]
Dividing both sides of the previous estimate by $\left|h\right|^{2}$,
using the properties of $\chi$, $\eta$ and $u_{\varepsilon}$, and
applying Young's inequality with exponents $\left(\frac{p}{p-2},\frac{p}{2}\right)$,
we get\begin{eqnarray}\label{eq:ee8}
\begin{split}
&\sup_{t\in(t_{0}-\ell^{2},t_{0})}\int_{B_{\ell}(x_{0})}\left|\Delta_{h}u_{\varepsilon}(x,t)\right|^{2}dx\,+\int_{Q_{\ell}(z_{0})}\left|\Delta_{h}H_{\frac{p}{2}}(Du_{\varepsilon})\right|^{2}dz+\,\varepsilon\,\int_{Q_{\ell}(z_{0})}\left|\Delta_{h}V_{p}(Du_{\varepsilon})\right|^{2}dz\\
&\leq\,\frac{c_{7}(p)}{(\rho-\ell)^{2}}\int_{Q_{\rho}(z_{0})}\left[\left(\left|H_{\frac{p}{2}}(Du_{\varepsilon}(x+he_{s},t))\right|^{\frac{p-2}{p}}+\left|H_{\frac{p}{2}}(Du_{\varepsilon})\right|^{\frac{p-2}{p}}\right)^{2}+1\right]\left|\Delta_{h}u_{\varepsilon}\right|^{2}dz\\
&\,\,\,\,\,\,\,+\,\,\frac{\varepsilon\,c_{7}(p)}{(\rho-\ell)^{2}}\int_{Q_{\rho}(z_{0})}\left(\left|Du_{\varepsilon}(x+he_{s},t)\right|^{\frac{p-2}{2}}+\left|Du_{\varepsilon}\right|^{\frac{p-2}{2}}\right)^{2}\left|\Delta_{h}u_{\varepsilon}\right|^{2}dz\\
&\,\,\,\,\,\,\,+\,\,c_{6}(n,p)\,\left(\int_{Q_{\gamma}(z_{0})}\left|Df_{\varepsilon}\right|^{\frac{np\,+\,4}{np\,+\,4\,-\,n}}dz\right)^{\frac{np\,+\,4\,-\,n}{np\,+\,4}}\left(\int_{Q_{\gamma}\left(z_{0}\right)}\left|Du_{\varepsilon}\right|^{p\,+\,\frac{4}{n}}dz\right)^{\frac{n}{np\,+\,4}}\\
&\leq\,\frac{c_{8}(p)}{(\rho-\ell)^{2}}\int_{Q_{\gamma}(z_{0})}\left(\left|Du_{\varepsilon}\right|^{p-2}+1\right)\left|\Delta_{h}u_{\varepsilon}\right|^{2}dz\\
&\,\,\,\,\,\,\,+\,\,c_{6}(n,p)\,\left(\int_{Q_{\gamma}(z_{0})}\left|Df_{\varepsilon}\right|^{\frac{np\,+\,4}{np\,+\,4\,-\,n}}dz\right)^{\frac{np\,+\,4\,-\,n}{np\,+\,4}}\left(\int_{Q_{\gamma}\left(z_{0}\right)}\left|Du_{\varepsilon}\right|^{p\,+\,\frac{4}{n}}dz\right)^{\frac{n}{np\,+\,4}}\\
&=\,\,\frac{c_{8}(p)}{(\rho-\ell)^{2}}\int_{Q_{\gamma}(z_{0})}\left|Du_{\varepsilon}\right|^{p-2}\left|\Delta_{h}u_{\varepsilon}\right|^{2}dz\,+\frac{c_{8}(p)}{(\rho-\ell)^{2}}\int_{Q_{\gamma}(z_{0})}\left|\Delta_{h}u_{\varepsilon}\right|^{2}dz\\
&\,\,\,\,\,\,\,+\,\,c_{6}(n,p)\,\left(\int_{Q_{\gamma}(z_{0})}\left|Df_{\varepsilon}\right|^{\frac{np\,+\,4}{np\,+\,4\,-\,n}}dz\right)^{\frac{np\,+\,4\,-\,n}{np\,+\,4}}\left(\int_{Q_{\gamma}\left(z_{0}\right)}\left|Du_{\varepsilon}\right|^{p\,+\,\frac{4}{n}}dz\right)^{\frac{n}{np\,+\,4}}\\
&\leq\,\frac{c_{9}(p)}{(\rho-\ell)^{2}}\int_{Q_{\gamma}}\left|Du_{\varepsilon}\right|^{p}dz\,+\,\frac{c_{10}(p)}{(\rho-\ell)^{2}}\int_{Q_{\gamma}}\left|\Delta_{h}u_{\varepsilon}\right|^{p}dz\,+\,\frac{c_{11}(n,p)}{(\rho-\ell)^{2}}\int_{Q_{R}}\left|Du_{\varepsilon}\right|^{2}dz\\
&\,\,\,\,\,\,\,+\,\,c_{6}(n,p)\,\left(\int_{Q_{\gamma}(z_{0})}\left|Df_{\varepsilon}\right|^{\frac{np\,+\,4}{np\,+\,4\,-\,n}}dz\right)^{\frac{np\,+\,4\,-\,n}{np\,+\,4}}\left(\int_{Q_{\gamma}\left(z_{0}\right)}\left|Du_{\varepsilon}\right|^{p\,+\,\frac{4}{n}}dz\right)^{\frac{n}{np\,+\,4}}\\
&\leq\,\,\frac{c_{12}(n,p)}{(\rho-\ell)^{2}}\,\left(\int_{Q_{R}(z_{0})}\left|Du_{\varepsilon}\right|^{p}dz\,+\int_{Q_{R}(z_{0})}\left|Du_{\varepsilon}\right|^{2}dz\right)\\
&\,\,\,\,\,\,\,+\,\,c_{12}(n,p)\,\left(\int_{Q_{R}(z_{0})}\left|Df_{\varepsilon}\right|^{\frac{np\,+\,4}{np\,+\,4\,-\,n}}dz\right)^{\frac{np\,+\,4\,-\,n}{np\,+\,4}}\left(\int_{Q_{\gamma}\left(z_{0}\right)}\left|Du_{\varepsilon}\right|^{p\,+\,\frac{4}{n}}dz\right)^{\frac{n}{np\,+\,4}},
\end{split}
\end{eqnarray}\\
where we have used Lemma \ref{lem:Giusti1}. Therefore, letting $h\rightarrow0$
in \eqref{eq:ee8}, we obtain\begin{equation}\label{eq:ee9}
\begin{split}
&\sup_{t\in(t_{0}-\ell^{2},t_{0})}\Vert Du_{\varepsilon}(\cdot,t)\Vert_{L^{2}(B_{\ell}(x_{0}))}^{2}+\int_{Q_{\ell}(z_{0})}\left|DH_{\frac{p}{2}}(Du_{\varepsilon})\right|^{2}dz\,+\,\varepsilon\,\int_{Q_{\ell}(z_{0})}\left|DV_{p}(Du_{\varepsilon})\right|^{2}dz\\
&\,\,\,\,\,\,\,\,\,\,\,\,\,\,\,\,\,\,\,\,\,\,\,\,\leq\,\,\frac{c_{12}(n,p)}{(\rho-\ell)^{2}}\,\left(\int_{Q_{R}(z_{0})}\left|Du_{\varepsilon}\right|^{p}dz\,+\int_{Q_{R}(z_{0})}\left|Du_{\varepsilon}\right|^{2}dz\right)\\
&\,\,\,\,\,\,\,\,\,\,\,\,\,\,\,\,\,\,\,\,\,\,\,\,\,\,\,\,\,\,\,+\,\,c_{12}(n,p)\,\left(\int_{Q_{R}(z_{0})}\left|Df_{\varepsilon}\right|^{\frac{np\,+\,4}{np\,+\,4\,-\,n}}dz\right)^{\frac{np\,+\,4\,-\,n}{np\,+\,4}}\left(\int_{Q_{\gamma}\left(z_{0}\right)}\left|Du_{\varepsilon}\right|^{p\,+\,\frac{4}{n}}dz\right)^{\frac{n}{np\,+\,4}}.
\end{split}
\end{equation}\\
Now we use \eqref{eq:ee3bis} in order to estimate the last integral
as follows\\
\begin{align*}
&\int_{Q_{\gamma}\left(z_{0}\right)}\left|Du_{\varepsilon}\right|^{p\,+\,\frac{4}{n}}dz\\
&=\int_{Q_{\gamma}\cap\left\{ \left|Du_{\varepsilon}\right|\,\geq\,\nu\right\} }\left(\left|Du_{\varepsilon}\right|-\nu+\nu\right)^{p\,+\,\frac{4}{n}}dz\,+\,\int_{Q_{\gamma}\cap\left\{ \left|Du_{\varepsilon}\right|\,<\,\nu\right\} }\left|Du_{\varepsilon}\right|^{p\,+\,\frac{4}{n}}dz\\
&\leq\,c_{13}(n,p)\int_{Q_{\gamma}}\left(\left|Du_{\varepsilon}\right|-\nu\right)_{+}^{p\,+\,\frac{4}{n}}dz\,+\,c_{13}(n,p)\,\nu^{p\,+\,\frac{4}{n}}\left|Q_{R_{0}}\right|\\
&\leq\,c_{14}(n,p)\left[\sup_{t\in(t_{0}-(\lambda r)^{2},t_{0})}\int_{B_{\lambda r}(x_{0})}\left|Du_{\varepsilon}(x,t)\right|^{2}dx\right]^{\frac{2}{n}}\cdot\int_{Q_{\lambda r}(z_{0})}\left(\left|DH_{\frac{p}{2}}(Du_{\varepsilon})\right|^{2}+\frac{\left|Du_{\varepsilon}\right|^{p}}{(\lambda r-\gamma)^{2}}\right)dz\\
&\,\,\,\,\,\,\,+\,c_{14}(n,p)\,\nu^{p\,+\,\frac{4}{n}}\left|Q_{R_{0}}\right|,
\end{align*}where $\left|Q_{R_{0}}\right|$ denotes the $(n+1)$-dimensional Lebesgue
measure of the parabolic cylinder $Q_{R_{0}}(z_{0})$. Plugging the
above estimate into \eqref{eq:ee9} and applying Young's inequality
with $\theta\in(0,1)$ and exponents $\left(\frac{np\,+\,4}{np\,+\,2\,-\,n},\frac{np\,+\,4}{2},\frac{np\,+\,4}{n}\right)$,
we get\\
\begin{align*}
&\sup_{t\in(t_{0}-\ell^{2},t_{0})}\Vert Du_{\varepsilon}(\cdot,t)\Vert_{L^{2}(B_{\ell}(x_{0}))}^{2}+\int_{Q_{\ell}(z_{0})}\left|DH_{\frac{p}{2}}(Du_{\varepsilon})\right|^{2}dz\\
&\leq\,\,\frac{c_{12}(n,p)}{(\rho-\ell)^{2}}\,\left(\Vert Du_{\varepsilon}\Vert_{L^{p}(Q_{R})}^{p}\,+\,\Vert Du_{\varepsilon}\Vert_{L^{2}(Q_{R})}^{2}\right)\,+\,c_{15}(n,p)\,\nu\,\Vert Df_{\varepsilon}\Vert_{L^{\frac{np+4}{np+4-n}}(Q_{R})}\left|Q_{R_{0}}\right|^{\frac{n}{np\,+\,4}}\\
&\,\,\,\,\,\,\,+\,c_{15}(n,p)\,\Vert Df_{\varepsilon}\Vert_{L^{\frac{np+4}{np+4-n}}(Q_{R})}\left[\sup_{t\in(t_{0}-(\lambda r)^{2},t_{0})}\Vert Du_{\varepsilon}(\cdot,t)\Vert_{L^{2}(B_{\lambda r}(x_{0}))}^{2}\right]^{\frac{2}{np\,+\,4}}\cdot\\
&\,\,\,\,\,\,\,\,\,\,\,\,\,\,\,\,\,\,\,\,\,\,\,\,\,\,\,\,\,\,\,\,\,\,\,\,\,\,\,\,\,\,\,\,\,\,\,\,\,\,\,\,\,\,\,\,\,\,\,\,\,\,\,\,\,\,\,\,\,\,\,\,\,\,\,\,\,\,\cdot\left[\int_{Q_{\lambda r}(z_{0})}\left(\left|DH_{\frac{p}{2}}(Du_{\varepsilon})\right|^{2}+\,\frac{\left|Du_{\varepsilon}\right|^{p}}{(\lambda r-\gamma)^{2}}\right)dz\right]^{\frac{n}{np\,+\,4}}\\
&\leq\,\,\frac{c_{12}(n,p)}{(\rho-\ell)^{2}}\,\left(\Vert Du_{\varepsilon}\Vert_{L^{p}(Q_{R})}^{p}\,+\,\Vert Du_{\varepsilon}\Vert_{L^{2}(Q_{R})}^{2}\right)\,+\,c_{16}(n,p,R_{0})\,\nu\,\Vert Df_{\varepsilon}\Vert_{L^{\frac{np+4}{np+4-n}}(Q_{R})}\\
&\,\,\,\,\,\,\,+\,\,c_{17}(n,p)\,\theta^{\frac{n\,+\,2}{n\,-np\,-2}}\,\Vert Df_{\varepsilon}\Vert_{L^{\frac{np+4}{np+4-n}}(Q_{R})}^{\frac{np+4}{np+2-n}}+\,\theta\,\left[\sup_{t\in(t_{0}-(\lambda r)^{2},t_{0})}\Vert Du_{\varepsilon}(\cdot,t)\Vert_{L^{2}(B_{\lambda r}(x_{0}))}^{2}\right]\\
&\,\,\,\,\,\,\,+\,\,\theta\int_{Q_{\lambda r}(z_{0})}\left|DH_{\frac{p}{2}}(Du_{\varepsilon})\right|^{2}dz\,+\,\frac{\theta}{(\lambda r-\gamma)^{2}}\,\Vert Du_{\varepsilon}\Vert_{L^{p}(Q_{R})}^{p}.
\end{align*}Now, if we choose $\theta=\frac{1}{2}$ and set
\[
\Psi(\xi)\,:=\sup_{t\in(t_{0}-\xi^{2},t_{0})}\Vert Du_{\varepsilon}(\cdot,t)\Vert_{L^{2}(B_{\xi}(x_{0}))}^{2}+\int_{Q_{\xi}(z_{0})}\left|DH_{\frac{p}{2}}(Du_{\varepsilon})\right|^{2}dz,
\]
then the previous estimate turns into\begin{equation}\label{eq:ee10}
\begin{split}
\Psi(r)&\leq\,\Psi(\ell)\,\leq\,\frac{1}{2}\,\Psi(\lambda r)\,+\,\frac{c_{12}(n,p)}{(\rho-\ell)^{2}}\,\left(\Vert Du_{\varepsilon}\Vert_{L^{p}(Q_{R})}^{p}\,+\,\Vert Du_{\varepsilon}\Vert_{L^{2}(Q_{R})}^{2}\right)\\
&\,\,\,\,\,\,\,+\,\frac{1}{2\,(\lambda r-\gamma)^{2}}\,\Vert Du_{\varepsilon}\Vert_{L^{p}(Q_{R})}^{p}\,+\,c_{18}(n,p,R_{0})\left(\nu\,\Vert Df_{\varepsilon}\Vert_{L^{\frac{np+4}{np+4-n}}(Q_{R})}+\,\Vert Df_{\varepsilon}\Vert_{L^{\frac{np+4}{np+4-n}}(Q_{R})}^{\frac{np+4}{np+2-n}}\right).
\end{split}
\end{equation}Since \eqref{eq:ee10} holds for all $\frac{R}{2}\leq r<\ell<\rho<\gamma<\lambda r<R$,
with $1<\lambda<2$, we can now choose 
\[
\ell=r+\frac{1}{4}(\lambda r-r),\,\,\,\,\,\,\,\,\,\rho=r+\frac{1}{2}(\lambda r-r)\,\,\,\,\,\,\,\,\,\mathrm{and}\,\,\,\,\,\,\,\,\,\gamma=r+\frac{3}{4}(\lambda r-r),
\]
thus obtaining\begin{align*}
\Psi(r)\,&\leq\,\frac{1}{2}\,\Psi(\lambda r)\,+\,\frac{c_{19}(n,p)}{r^{2}(\lambda-1)^{2}}\,\left(\Vert Du_{\varepsilon}\Vert_{L^{p}(Q_{R})}^{p}\,+\,\Vert Du_{\varepsilon}\Vert_{L^{2}(Q_{R})}^{2}\right)\\
&\,\,\,\,\,\,\,+\,\,c_{18}(n,p,R_{0})\left(\nu\,\Vert Df_{\varepsilon}\Vert_{L^{\frac{np+4}{np+4-n}}(Q_{R})}+\,\Vert Df_{\varepsilon}\Vert_{L^{\frac{np+4}{np+4-n}}(Q_{R})}^{\frac{np+4}{np+2-n}}\right),
\end{align*}which holds for all $\frac{R}{2}\leq r<\lambda r<R$, with $1<\lambda<2$.
Therefore, the use of Lemma \ref{lem:giusti} with $r_{0}=\frac{R}{2}$
and $r_{1}=R$ yields\\
\begin{align*}
\Psi(R/2)\,&\leq\,\frac{c_{20}(n,p)}{R^{2}}\,\left(\Vert Du_{\varepsilon}\Vert_{L^{p}(Q_{R})}^{p}\,+\,\Vert Du_{\varepsilon}\Vert_{L^{2}(Q_{R})}^{2}\right)\\
&\,\,\,\,\,\,\,+\,c_{21}(n,p,R_{0})\left(\nu\,\Vert Df_{\varepsilon}\Vert_{L^{\frac{np+4}{np+4-n}}(Q_{R})}+\,\Vert Df_{\varepsilon}\Vert_{L^{\frac{np+4}{np+4-n}}(Q_{R})}^{\frac{np+4}{np+2-n}}\right),
\end{align*}that is,\begin{equation}\label{eq:essenz1}
\begin{split}
&\sup_{t\in(t_{0}-R^{2}/4,t_{0})}\Vert Du_{\varepsilon}(\cdot,t)\Vert_{L^{2}(B_{R/2}(x_{0}))}^{2}+\int_{Q_{R/2}(z_{0})}\left|DH_{\frac{p}{2}}(Du_{\varepsilon})\right|^{2}dz\\
&\,\,\,\,\,\,\,\,\,\,\,\,\,\,\,\,\,\,\,\,\,\,\,\,\,\,\,\,\,\,\leq\,\frac{c_{20}(n,p)}{R^{2}}\left(\Vert Du_{\varepsilon}\Vert_{L^{p}(Q_{R})}^{p}\,+\,\Vert Du_{\varepsilon}\Vert_{L^{2}(Q_{R})}^{2}\right)\\
&\,\,\,\,\,\,\,\,\,\,\,\,\,\,\,\,\,\,\,\,\,\,\,\,\,\,\,\,\,\,\,\,\,\,\,\,\,+\,c_{21}(n,p,R_{0})\left(\nu\,\Vert Df_{\varepsilon}\Vert_{L^{\frac{np+4}{np+4-n}}(Q_{R})}+\,\Vert Df_{\varepsilon}\Vert_{L^{\frac{np+4}{np+4-n}}(Q_{R})}^{\frac{np+4}{np+2-n}}\right),
\end{split}
\end{equation}\\
which is the \textit{a priori} estimate we were looking for. \\
$\hspace*{1em}$Our next aim is to prove that the norms $\Vert Du_{\varepsilon}\Vert_{L^{p}(Q_{R_{0}})}$
are all bounded by a constant independent of $\varepsilon$ (and therefore
the same will be true for the $L^{2}$-norms, since $p>2$).\textbf{}\\
\textbf{}\\
\textbf{}\\
\textbf{Step 2: the uniform boundedness of the norms $\boldsymbol{\Vert Du_{\varepsilon}\Vert_{L^{p}(Q_{R_{0}})}}$.}\\
\\
In order to have an (uniform in $\varepsilon$) energy estimate for
$\vert Du_{\varepsilon}\vert^{p}$, we now proceed by testing equations
(\ref{eq:1}) and $(\ref{eq:Diri})_{1}$ for $u$ and $u_{\varepsilon}$,
respectively, with the map $\varphi=g(t)(u_{\varepsilon}-u)$, where
$g\in W^{1,\infty}(\mathbb{R})$ is chosen such that\vspace{0.3cm}
\[
g(t)=\begin{cases}
\begin{array}{c}
1\\
-\,\frac{1}{\omega}(t-t_{2}-\omega)\\
0
\end{array} & \begin{array}{c}
\mathrm{if}\,\,t\leq t_{2},\quad\quad\quad\;\,\,\,\,\\
\mathrm{if}\,\,t_{2}<t<t_{2}+\omega,\\
\mathrm{if}\,\,t\geq t_{2}+\omega,\quad\quad
\end{array}\end{cases}
\]
\\
with $t_{0}-R_{0}^{2}<t_{2}<t_{2}\,+\,\omega<t_{0}$, and then letting
$\omega\rightarrow0$. We observe that at this stage it is important
that $u_{\varepsilon}$ and $u$ agree on the parabolic boundary $\partial_{\mathrm{par}}Q_{R_{0}}(z_{0})$.
We also note that the following computations are somewhat formal concerning
the use of the time derivative, but they can easily be made rigorous,
for example by the use of Steklov averages. We skip this, since it
is a standard procedure. With the previous choice of $\varphi$, for
almost every $t_{2}\in(t_{0}\,-\,R_{0}^{2},t_{0})$ we find\vspace{0.3cm}
\[
\frac{1}{2}\int_{B_{R_{0}}(x_{0})}\left|u_{\varepsilon}(x,t_{2})-u(x,t_{2})\right|^{2}dx\,+\int_{Q_{R_{0},t_{2}}}\langle H_{p-1}(Du_{\varepsilon})-H_{p-1}(Du),Du_{\varepsilon}-Du\rangle\,dz
\]
\begin{equation}
+\,\,\varepsilon\int_{Q_{R_{0},t_{2}}}\langle\vert Du_{\varepsilon}\vert^{p-2}Du_{\varepsilon},Du_{\varepsilon}-Du\rangle\,dz\,=\int_{Q_{R_{0},t_{2}}}(f-f_{\varepsilon})(u_{\varepsilon}-u)\,dz,\label{eq:ee11}
\end{equation}
\\
where we have used the abbreviation $Q_{R_{0},t_{2}}=B_{R_{0}}(x_{0})\times(t_{0}-R_{0}^{2},t_{2})$.
Using Lemma \ref{lem:Brasco}, the Cauchy-Schwarz inequality as well
as Young's inequality with $\mu>0$ and exponents $(p,p')$, from
(\ref{eq:ee11}) we infer\\
\\
\begin{equation}\label{eq:ee12}
\begin{split}
&\frac{1}{2}\sup_{t\in(t_{0}-R_{0}^{2},t_{0})}\Vert u_{\varepsilon}(\cdot,t)-u(\cdot,t)\Vert_{L^{2}(B_{R_{0}}(x_{0}))}^{2}\,+\,\frac{4}{p^{2}}\int_{Q_{R_{0}}(z_{0})}\left|H_{\frac{p}{2}}(Du_{\varepsilon})-H_{\frac{p}{2}}(Du)\right|^{2}dz\\
&+\,\,\varepsilon\int_{Q_{R_{0}}(z_{0})}\left|Du_{\varepsilon}\right|^{p}dz\\
&\,\,\,\,\,\,\,\leq\,\int_{Q_{R_{0}}(z_{0})}\left|f-f_{\varepsilon}\right|\left|u_{\varepsilon}-u\right|\,dz\,+\,\varepsilon\int_{Q_{R_{0}}(z_{0})}\left|Du_{\varepsilon}\right|^{p-1}\left|Du\right|\,dz\\
&\,\,\,\,\,\,\,\leq\,\int_{Q_{R_{0}}}\left|f-f_{\varepsilon}\right|\left|u_{\varepsilon}-u\right|\,dz\,+\,\frac{\varepsilon}{p\,\mu^{p/2}}\int_{Q_{R_{0}}}\left|Du\right|^{p}dz\,+\,\frac{\varepsilon}{p'}\,\mu^{\frac{p'}{2}}\int_{Q_{R_{0}}}\left|Du_{\varepsilon}\right|^{p}dz.
\end{split}
\end{equation}\\
Choosing $\mu=\left(\frac{p'}{2}\right)^{\frac{2}{p'}}$ and reabsorbing
the last integral in the right-hand side of \eqref{eq:ee12} by the
left-hand side, we arrive at\\
\\
\begin{equation}\label{eq:ee13}
\begin{split}
&\sup_{t\in(t_{0}-R_{0}^{2},t_{0})}\Vert u_{\varepsilon}(\cdot,t)-u(\cdot,t)\Vert_{L^{2}(B_{R_{0}}(x_{0}))}^{2}\,+\int_{Q_{R_{0}}}\left|H_{\frac{p}{2}}(Du_{\varepsilon})-H_{\frac{p}{2}}(Du)\right|^{2}dz\,+\,\varepsilon\int_{Q_{R_{0}}}\left|Du_{\varepsilon}\right|^{p}dz\\
&\,\,\,\,\,\,\,\leq\,\varepsilon\,c_{1}(p)\int_{Q_{R_{0}}(z_{0})}\left|Du\right|^{p}dz\,+\,c_{1}(p)\int_{Q_{R_{0}}(z_{0})}\left|f-f_{\varepsilon}\right|\left|u_{\varepsilon}-u\right|\,dz.
\end{split}
\end{equation}Now we set\vspace{0.5cm}
\[
A_{5}:=\int_{Q_{R_{0}}(z_{0})}\left|f-f_{\varepsilon}\right|\left|u_{\varepsilon}-u\right|\,dz
\]
and observe that\vspace{0.5cm}
\begin{equation}
\left(p+\frac{4}{n}\right)'=\frac{np\,+\,4}{np\,+4-n}\,\geq\,\left(p+\frac{2p}{n}\right)'=\frac{np\,+\,2p}{np\,+\,2p\,-n}\,\,\,\,\,\,\,\,\,\,\,\mathrm{for\,\,every\,\,}p\geq2.\label{eq:expon}
\end{equation}
\\
\\
Therefore, we can apply Hölder's inequality with exponents $\left(\frac{np\,+\,2p}{np\,+\,2p\,-\,n},p+\frac{2p}{n}\right)$
and Lemma \ref{lem:inter} with $v=u_{\varepsilon}-u$, $r=R_{0}$
and $q=2$ to estimate $A_{5}$ as follows:\\
\\
\begin{equation}\label{eq:ee14}
\begin{split}
A_{5}\,&\leq\,\left(\int_{Q_{R_{0}}(z_{0})}\left|f-f_{\varepsilon}\right|^{\frac{np\,+\,2p}{np\,+\,2p\,-n}}dz\right)^{\frac{np\,+\,2p\,-n}{np\,+\,2p}}\left(\int_{Q_{R_{0}}(z_{0})}\left|u_{\varepsilon}-u\right|^{p\,+\,\frac{2p}{n}}dz\right)^{\frac{n}{np\,+\,2p}}\\
&\leq\,c_{2}(n,p)\,\Vert f-f_{\varepsilon}\Vert_{L^{\frac{np\,+\,2p}{np\,+\,2p\,-n}}(Q_{R_{0}})}\left(\sup_{t\in(t_{0}-R_{0}^{2},t_{0})}\Vert u_{\varepsilon}(\cdot,t)-u(\cdot,t)\Vert_{L^{2}(B_{R_{0}}(x_{0}))}^{2}\right)^{\frac{1}{n\,+\,2}}\cdot\\
&\,\,\,\,\,\,\,\,\,\,\,\,\,\,\,\,\,\,\,\,\,\,\,\,\,\,\,\,\,\,\,\,\,\,\,\,\,\,\,\,\,\,\,\,\,\,\,\,\,\,\,\,\,\,\,\,\,\,\,\,\,\,\,\,\,\,\,\,\,\,\,\,\,\,\,\,\,\,\,\,\,\cdot\left(\int_{Q_{R_{0}}(z_{0})}\left|Du_{\varepsilon}-Du\right|^{p}dz\right)^{\frac{n}{np\,+\,2p}}.
\end{split}
\end{equation}\\
Combining estimates \eqref{eq:ee13} and \eqref{eq:ee14}, recalling
that $0<\varepsilon\leq1$ and applying Young's inequality with $\beta>0$
and exponents $\left(\frac{np\,+\,2p}{np\,+\,p\,-\,n},n+2,\frac{np\,+\,2p}{n}\right)$,
we obtain\\
\begin{equation}\label{eq:ee15}
\begin{split}
&\sup_{t\in(t_{0}-R_{0}^{2},t_{0})}\Vert u_{\varepsilon}(\cdot,t)-u(\cdot,t)\Vert_{L^{2}(B_{R_{0}}(x_{0}))}^{2}\,+\int_{Q_{R_{0}}}\left|H_{\frac{p}{2}}(Du_{\varepsilon})-H_{\frac{p}{2}}(Du)\right|^{2}dz\,+\,\varepsilon\int_{Q_{R_{0}}}\left|Du_{\varepsilon}\right|^{p}dz\\
&\,\,\,\,\,\,\,\leq\,c_{1}(p)\int_{Q_{R_{0}}}\left|Du\right|^{p}dz\,+\,c_{3}(n,p)\,\beta^{\frac{n\,+\,p}{n-np-p}}\,\Vert f-f_{\varepsilon}\Vert_{L^{\frac{np\,+\,2p}{np\,+\,2p\,-n}}(Q_{R_{0}})}^{\frac{np\,+\,2p}{np\,+\,p\,-n}}\\
&\,\,\,\,\,\,\,\,\,\,\,\,\,\,+\,\,\beta\left[\sup_{t\in(t_{0}-R_{0}^{2},t_{0})}\Vert u_{\varepsilon}(\cdot,t)-u(\cdot,t)\Vert_{L^{2}(B_{R_{0}}(x_{0}))}^{2}\right]+\,\beta\int_{Q_{R_{0}}}\left|Du_{\varepsilon}-Du\right|^{p}dz\\
&\,\,\,\,\,\,\,\leq\,\left(c_{1}(p)+2^{p-1}\beta\right)\int_{Q_{R_{0}}}\left|Du\right|^{p}dz\,+\,c_{3}(n,p)\,\beta^{\frac{n\,+\,p}{n-np-p}}\,\Vert f-f_{\varepsilon}\Vert_{L^{\frac{np\,+\,2p}{np\,+\,2p\,-n}}(Q_{R_{0}})}^{\frac{np\,+\,2p}{np\,+\,p\,-n}}\\
&\,\,\,\,\,\,\,\,\,\,\,\,\,\,+\,\,\beta\left[\sup_{t\in(t_{0}-R_{0}^{2},t_{0})}\Vert u_{\varepsilon}(\cdot,t)-u(\cdot,t)\Vert_{L^{2}(B_{R_{0}}(x_{0}))}^{2}\right]+\,2^{p-1}\beta\int_{Q_{R_{0}}}\left|Du_{\varepsilon}\right|^{p}dz.
\end{split}
\end{equation}Now, let us notice that

\noindent \begin{equation}\label{eq:ee16}
\begin{split}
&\int_{Q_{R_{0}}(z_{0})}\left|Du_{\varepsilon}\right|^{p}dz\\
&\,\,\,\,\,\,\,=\int_{Q_{R_{0}}\cap\left\{ \left|Du_{\varepsilon}\right|\,\geq\,\nu\right\} }\left(\left|Du_{\varepsilon}\right|-\nu+\nu\right)^{p}dz\,+\,\int_{Q_{R_{0}}\cap\left\{ \left|Du_{\varepsilon}\right|\,<\,\nu\right\} }\left|Du_{\varepsilon}\right|^{p}dz\\
&\,\,\,\,\,\,\,\leq\,2^{p-1}\int_{Q_{R_{0}}}\left[\left(\left|Du_{\varepsilon}\right|-\nu\right)_{+}^{p}+\nu^{p}\right]\,dz\,+\,\nu^{p}\left|Q_{R_{0}}\right|\\
&\,\,\,\,\,\,\,=\,2^{p-1}\int_{Q_{R_{0}}}\left(\left|H_{\frac{p}{2}}(Du_{\varepsilon})-H_{\frac{p}{2}}(Du)+H_{\frac{p}{2}}(Du)\right|^{2}+\nu^{p}\right)\,dz\,+\,\nu^{p}\left|Q_{R_{0}}\right|\\
&\,\,\,\,\,\,\,\leq\,2^{p}\int_{Q_{R_{0}}}\left(\left|H_{\frac{p}{2}}(Du_{\varepsilon})-H_{\frac{p}{2}}(Du)\right|^{2}+\left|H_{\frac{p}{2}}(Du)\right|^{2}+\nu^{p}\right)\,dz\,+\,\nu^{p}\left|Q_{R_{0}}\right|\\
&\,\,\,\,\,\,\,\leq\,2^{p}\int_{Q_{R_{0}}}\left|H_{\frac{p}{2}}(Du_{\varepsilon})-H_{\frac{p}{2}}(Du)\right|^{2}dz\,+\,2^{p+1}\int_{Q_{R_{0}}}\left(\left|Du\right|^{p}+\nu^{p}\right)\,dz.
\end{split}
\end{equation}\\
\\
Multiplying all sides of \eqref{eq:ee15} by $2^{p}$ and then adding
the resulting expression and \eqref{eq:ee16} side by side, we get\\
\\
\begin{equation}\label{eq:ee17}
\begin{split}
&(1+2^{p}\varepsilon)\int_{Q_{R_{0}}}\left|Du_{\varepsilon}\right|^{p}dz\,+\,2^{p}\sup_{t\in(t_{0}-R_{0}^{2},t_{0})}\Vert u_{\varepsilon}(\cdot,t)-u(\cdot,t)\Vert_{L^{2}(B_{R_{0}}(x_{0}))}^{2}\\
&\,\,\,\,\,\,\,\leq\,\left(c_{4}(p)+2^{2p-1}\beta\right)\int_{Q_{R_{0}}}\left(\left|Du\right|^{p}+\nu^{p}\right)dz\,+\,c_{5}(n,p)\,\beta^{\frac{n\,+\,p}{n-np-p}}\,\Vert f-f_{\varepsilon}\Vert_{L^{\frac{np\,+\,2p}{np\,+\,2p\,-n}}(Q_{R_{0}})}^{\frac{np\,+\,2p}{np\,+\,p\,-n}}\\
&\,\,\,\,\,\,\,\,\,\,\,\,\,\,+\,\,2^{p}\beta\left[\sup_{t\in(t_{0}-R_{0}^{2},t_{0})}\Vert u_{\varepsilon}(\cdot,t)-u(\cdot,t)\Vert_{L^{2}(B_{R_{0}}(x_{0}))}^{2}\right]+\,2^{2p-1}\beta\int_{Q_{R_{0}}}\left|Du_{\varepsilon}\right|^{p}dz.
\end{split}
\end{equation}\\
Now, choosing $\beta=1/2^{2p}$ and reabsorbing the last two terms
in the right-hand side of \eqref{eq:ee17} by the left-hand side,
we obtain\begin{align*}
&\int_{Q_{R_{0}}}\left|Du_{\varepsilon}\right|^{p}dz\,+\sup_{t\in(t_{0}-R_{0}^{2},t_{0})}\Vert u_{\varepsilon}(\cdot,t)-u(\cdot,t)\Vert_{L^{2}(B_{R_{0}}(x_{0}))}^{2}\\
&\,\,\,\,\,\,\,\leq\,c_{6}(p)\int_{Q_{R_{0}}}\left(\left|Du\right|^{p}+\nu^{p}\right)\,dz\,+\,c_{7}(n,p)\,\Vert f-f_{\varepsilon}\Vert_{L^{\frac{np\,+\,2p}{np\,+\,2p\,-n}}(Q_{R_{0}})}^{\frac{np\,+\,2p}{np\,+\,p\,-n}}. 
\end{align*}Moreover, since $f\in L^{\frac{np\,+\,4}{np\,+\,4\,-\,n}}(\Omega_{T})$,
recalling (\ref{eq:molli}), by virtue of (\ref{eq:expon}) we have\\
\begin{equation}
f_{\varepsilon}\rightarrow f\,\,\,\,\,\,\,\,\,\,\,\,\mathrm{strongly\,\,in\,\,}L^{\frac{np\,+\,2p}{np\,+\,2p\,-\,n}}(Q_{R_{0}}),\,\,\,\,\,\mathrm{as}\,\,\varepsilon\rightarrow0.\label{eq:strong}
\end{equation}
\\
Hence, there exists a finite positive number $M$, independent of
$\varepsilon$, such that\\
\[
\Vert f-f_{\varepsilon}\Vert_{L^{\frac{np\,+\,2p}{np\,+\,2p\,-n}}(Q_{R_{0}})}\leq\,M\,\,\,\,\,\,\,\,\,\,\,\,\,\,\mathrm{for\,\,all}\,\,\varepsilon\in(0,1],
\]
\\
from which we can finally deduce\medskip{}
\begin{equation}\label{eq:ee18}
\begin{split}
&\Vert Du_{\varepsilon}\Vert_{L^{p}(Q_{R_{0}})}^{p}\,+\sup_{t\in(t_{0}-R_{0}^{2},t_{0})}\Vert u_{\varepsilon}(\cdot,t)-u(\cdot,t)\Vert_{L^{2}(B_{R_{0}}(x_{0}))}^{2}\\
&\,\,\,\,\,\,\,\leq\,c_{8}(n,p,R_{0})\,\left(\Vert Du\Vert_{L^{p}(Q_{R_{0}})}^{p}+\nu^{p}+\Vert f-f_{\varepsilon}\Vert_{L^{\frac{np\,+\,2p}{np\,+\,2p\,-n}}(Q_{R_{0}})}^{\frac{np\,+\,2p}{np\,+\,p\,-n}}\right)\leq K\,\,\,\,\,\,\,\,\,\,\,\,\,\mathrm{for\,\,every}\,\,\varepsilon\in(0,1],
\end{split}
\end{equation}for some positive constant $K$ depending on $n$, $p$, $\nu$, $M$,
$R_{0}$ and $\Vert Du\Vert_{L^{p}(Q_{R_{0}})}$, but not on $\varepsilon$.\\
Joining \eqref{eq:ee13}, \eqref{eq:ee14} and \eqref{eq:ee18}, and
applying Minkowski's inequality, we get the following\textbf{ }\textbf{\textit{comparison
estimate}}\begin{align*}
&\sup_{t\in(t_{0}-R_{0}^{2},t_{0})}\Vert u_{\varepsilon}(\cdot,t)-u(\cdot,t)\Vert_{L^{2}(B_{R_{0}}(x_{0}))}^{2}\,+\int_{Q_{R_{0}}}\left|H_{\frac{p}{2}}(Du_{\varepsilon})-H_{\frac{p}{2}}(Du)\right|^{2}dz\\
&\,\,\,\,\,\,\,\,\,\,\,\,\,\,\leq\,\varepsilon\,c_{1}(p)\,\Vert Du\Vert_{L^{p}(Q_{R_{0}})}^{p}\\
&\,\,\,\,\,\,\,\,\,\,\,\,\,\,\,\,\,\,\,\,\,+\,c_{9}(n,p)\,\Vert f-f_{\varepsilon}\Vert_{L^{\frac{np\,+\,2p}{np\,+\,2p\,-n}}(Q_{R_{0}})}\left(\sup_{t\in(t_{0}-R_{0}^{2},t_{0})}\Vert u_{\varepsilon}(\cdot,t)-u(\cdot,t)\Vert_{L^{2}(B_{R_{0}}(x_{0}))}^{2}\right)^{\frac{1}{n\,+\,2}}\cdot\\
&\,\,\,\,\,\,\,\,\,\,\,\,\,\,\,\,\,\,\,\,\,\,\,\,\,\,\,\,\,\,\,\,\,\,\,\,\,\,\,\,\,\,\,\,\,\,\,\,\,\,\,\,\,\,\,\,\,\,\,\,\,\,\,\,\,\,\,\,\,\,\,\,\,\,\,\,\,\,\,\,\,\,\,\,\,\,\,\,\,\,\,\,\,\,\,\,\,\,\,\,\,\,\cdot\left(\int_{Q_{R_{0}}}\left|Du_{\varepsilon}-Du\right|^{p}dz\right)^{\frac{n}{np\,+\,2p}}
\end{align*}\begin{equation}\label{eq:comparison}
\begin{split}
&\,\,\,\,\,\,\,\,\,\,\,\,\,\,\leq\,\varepsilon\,c_{1}(p)\,\Vert Du\Vert_{L^{p}(Q_{R_{0}})}^{p}\\
&\,\,\,\,\,\,\,\,\,\,\,\,\,\,\,\,\,\,\,\,\,+\,c_{10}(n,p,R_{0})\,\Vert f-f_{\varepsilon}\Vert_{L^{\frac{np\,+\,2p}{np\,+\,2p\,-n}}(Q_{R_{0}})}\left(\Vert Du\Vert_{L^{p}(Q_{R_{0}})}^{p}+\nu^{p}+\Vert f-f_{\varepsilon}\Vert_{L^{\frac{np\,+\,2p}{np\,+\,2p\,-n}}(Q_{R_{0}})}^{\frac{np\,+\,2p}{np\,+\,p\,-n}}\right)^{^{\frac{1}{n+2}}}\cdot\\
&\,\,\,\,\,\,\,\,\,\,\,\,\,\,\,\,\,\,\,\,\,\,\,\,\,\,\,\,\,\,\,\,\,\,\,\,\,\,\,\,\,\,\,\,\,\,\,\,\,\,\,\,\,\,\,\,\,\,\,\,\,\,\,\,\,\,\,\,\,\,\,\,\,\,\,\,\,\,\,\,\,\,\,\,\,\,\,\,\,\,\,\,\,\,\,\,\,\,\,\,\,\,\,\,\,\,\,\,\,\,\,\,\,\,\cdot\left(\Vert Du_{\varepsilon}\Vert_{L^{p}(Q_{R_{0}})}+\Vert Du\Vert_{L^{p}(Q_{R_{0}})}\right)^{\frac{n}{n\,+\,2}}\\
&\,\,\,\,\,\,\,\,\,\,\,\,\,\,\leq\,\varepsilon\,c_{1}(p)\,\Vert Du\Vert_{L^{p}(Q_{R_{0}})}^{p}\\
&\,\,\,\,\,\,\,\,\,\,\,\,\,\,\,\,\,\,\,\,\,+\,c_{11}(n,p,R_{0})\,\Vert f-f_{\varepsilon}\Vert_{L^{\frac{np\,+\,2p}{np\,+\,2p\,-n}}(Q_{R_{0}})}\left(\Vert Du\Vert_{L^{p}(Q_{R_{0}})}^{p}+\nu^{p}+\Vert f-f_{\varepsilon}\Vert_{L^{\frac{np\,+\,2p}{np\,+\,2p\,-n}}(Q_{R_{0}})}^{\frac{np\,+\,2p}{np\,+\,p\,-n}}\right)^{^{\frac{1}{n+2}}}\cdot\\
&\,\,\,\,\,\,\,\,\,\,\,\,\,\,\,\,\,\,\,\,\,\,\,\,\,\,\,\,\,\,\,\,\,\,\,\,\,\,\,\,\,\,\,\,\,\,\,\,\,\,\,\,\,\,\,\,\,\,\,\,\,\,\,\,\,\,\,\,\,\,\,\,\,\,\,\,\,\,\,\,\,\,\,\,\,\,\,\,\,\,\,\,\,\,\,\,\,\,\,\,\,\,\,\,\,\,\,\,\,\,\,\,\,\cdot\left(\Vert Du\Vert_{L^{p}(Q_{R_{0}})}+\nu+\Vert f-f_{\varepsilon}\Vert_{L^{\frac{np\,+\,2p}{np\,+\,2p\,-n}}(Q_{R_{0}})}^{\frac{n+2}{np\,+\,p\,-n}}\right)^{\frac{n}{n\,+\,2}},
\end{split}
\end{equation}where, in the last line, we have used estimate \eqref{eq:ee18} again.\\
\\
\textbf{Step 3: the conclusion.}\\
\\
Now we shall fix $\varrho\in(0,R)$ and consider the finite difference
operator $\tau_{s,h}$, defined in Section \ref{subsec:DiffOpe},
for increments $h\in\mathbb{R}\setminus\left\{ 0\right\} $ such that
$\left|h\right|$ is suitably small.\\
In order to obtain an estimate for the finite difference $\tau_{s,h}H_{\frac{p}{2}}(Du)$,
we use the following comparison argument:\begin{align*}
&\int_{Q_{\varrho/2}(z_{0})}\left|\tau_{s,h}H_{\frac{p}{2}}(Du)\right|^{2}dz\,\leq\,4\int_{Q_{\varrho/2}(z_{0})}\left|\tau_{s,h}H_{\frac{p}{2}}(Du_{\varepsilon})\right|^{2}dx\,dt\\
&\,\,\,\,\,\,\,+\,\,4\int_{Q_{\varrho/2}(z_{0})}\left|H_{\frac{p}{2}}(Du_{\varepsilon}(x+he_{s},t))-H_{\frac{p}{2}}(Du(x+he_{s},t))\right|^{2}dx\,dt\\
&\,\,\,\,\,\,\,+\,\,4\int_{Q_{\varrho/2}(z_{0})}\left|H_{\frac{p}{2}}(Du_{\varepsilon})-H_{\frac{p}{2}}(Du)\right|^{2}dx\,dt.
\end{align*} Combining the above inequality with estimates \eqref{eq:essenz1},
\eqref{eq:ee18} and \eqref{eq:comparison}, for every $s\in\left\{ 1,\ldots,n\right\} $
we get\begin{align*}
&\int_{Q_{\varrho/2}(z_{0})}\left|\tau_{s,h}H_{\frac{p}{2}}(Du)\right|^{2}dz\\
&\,\,\,\,\,\,\,\leq\,4\int_{Q_{R/2}(z_{0})}\left|\tau_{s,h}H_{\frac{p}{2}}(Du_{\varepsilon})\right|^{2}dz\,+\,8\int_{Q_{R}(z_{0})}\left|H_{\frac{p}{2}}(Du_{\varepsilon})-H_{\frac{p}{2}}(Du)\right|^{2}dz\\
&\,\,\,\,\,\,\,\leq\,\frac{c_{1}(n,p)}{R^{2}}\left|h\right|^{2}\left(\Vert Du_{\varepsilon}\Vert_{L^{p}(Q_{R})}^{p}\,+\,\Vert Du_{\varepsilon}\Vert_{L^{2}(Q_{R})}^{2}\right)+\,\varepsilon\,c_{2}(p)\,\Vert Du\Vert_{L^{p}(Q_{R_{0}})}^{p}\\
&\,\,\,\,\,\,\,\,\,\,\,\,\,\,+\,c_{3}(n,p,R_{0})\left|h\right|^{2}\left(\nu\,\Vert Df_{\varepsilon}\Vert_{L^{\frac{np+4}{np+4-n}}(Q_{R})}+\,\Vert Df_{\varepsilon}\Vert_{L^{\frac{np+4}{np+4-n}}(Q_{R})}^{\frac{np+4}{np+2-n}}\right)\\
&\,\,\,\,\,\,\,\,\,\,\,\,\,\,+\,c_{3}(n,p,R_{0})\Vert f-f_{\varepsilon}\Vert_{L^{\frac{np\,+\,2p}{np\,+\,2p\,-n}}(Q_{R_{0}})}\left(\Vert Du\Vert_{L^{p}(Q_{R_{0}})}^{p}+\nu^{p}+\Vert f-f_{\varepsilon}\Vert_{L^{\frac{np\,+\,2p}{np\,+\,2p\,-n}}(Q_{R_{0}})}^{\frac{np\,+\,2p}{np\,+\,p\,-n}}\right)^{^{\frac{1}{n+2}}}\\
&\,\,\,\,\,\,\,\,\,\,\,\,\,\,\,\,\,\,\,\,\,\,\,\,\,\,\,\,\,\,\,\,\,\,\,\,\,\,\,\,\,\,\,\,\,\,\,\,\,\,\,\,\,\,\,\,\,\,\,\,\,\,\,\,\,\,\,\,\,\,\,\,\,\,\,\,\,\,\,\,\,\,\,\,\,\,\,\,\,\,\,\,\,\,\,\,\,\,\,\,\,\,\,\,\cdot\left(\Vert Du\Vert_{L^{p}(Q_{R_{0}})}+\nu+\Vert f-f_{\varepsilon}\Vert_{L^{\frac{np\,+\,2p}{np\,+\,2p\,-n}}(Q_{R_{0}})}^{\frac{n+2}{np\,+\,p\,-n}}\right)^{\frac{n}{n\,+\,2}}
\end{align*}\begin{equation}\label{eq:ee19}
\begin{split}
&\,\,\,\,\,\,\,\leq\,\frac{c_{4}(n,p,R_{0})}{R^{2}}\left|h\right|^{2}\left(\Vert Du\Vert_{L^{p}(Q_{R_{0}})}^{p}+\Vert Du\Vert_{L^{p}(Q_{R_{0}})}^{2}+\nu^{p}+\nu^{2}\right)\\
&\,\,\,\,\,\,\,\,\,\,\,\,\,\,+\,\frac{c_{4}(n,p,R_{0})}{R^{2}}\left|h\right|^{2}\left(\Vert f-f_{\varepsilon}\Vert_{L^{\frac{np\,+\,2p}{np\,+\,2p\,-n}}(Q_{R_{0}})}^{\frac{np\,+\,2p}{np\,+\,p\,-n}}+\Vert f-f_{\varepsilon}\Vert_{L^{\frac{np\,+\,2p}{np\,+\,2p\,-n}}(Q_{R_{0}})}^{\frac{2n\,+\,4}{np\,+\,p\,-n}}\right)\\
&\,\,\,\,\,\,\,\,\,\,\,\,\,\,+\,\varepsilon\,c_{2}(p)\,\Vert Du\Vert_{L^{p}(Q_{R_{0}})}^{p}+\,c_{4}(n,p,R_{0})\left|h\right|^{2}\left(\nu\,\Vert Df_{\varepsilon}\Vert_{L^{\frac{np+4}{np+4-n}}(Q_{R})}+\,\Vert Df_{\varepsilon}\Vert_{L^{\frac{np+4}{np+4-n}}(Q_{R})}^{\frac{np+4}{np+2-n}}\right)\\
&\,\,\,\,\,\,\,\,\,\,\,\,\,\,+\,c_{4}(n,p,R_{0})\Vert f-f_{\varepsilon}\Vert_{L^{\frac{np\,+\,2p}{np\,+\,2p\,-n}}(Q_{R_{0}})}\left(\Vert Du\Vert_{L^{p}(Q_{R_{0}})}^{p}+\nu^{p}+\Vert f-f_{\varepsilon}\Vert_{L^{\frac{np\,+\,2p}{np\,+\,2p\,-n}}(Q_{R_{0}})}^{\frac{np\,+\,2p}{np\,+\,p\,-n}}\right)^{^{\frac{1}{n+2}}}\\
&\,\,\,\,\,\,\,\,\,\,\,\,\,\,\,\,\,\,\,\,\,\,\,\,\,\,\,\,\,\,\,\,\,\,\,\,\,\,\,\,\,\,\,\,\,\,\,\,\,\,\,\,\,\,\,\,\,\,\,\,\,\,\,\,\,\,\,\,\,\,\,\,\,\,\,\,\,\,\,\,\,\,\,\,\,\,\,\,\,\,\,\,\,\,\,\,\,\,\,\,\,\,\,\,\cdot\left(\Vert Du\Vert_{L^{p}(Q_{R_{0}})}+\nu+\Vert f-f_{\varepsilon}\Vert_{L^{\frac{np\,+\,2p}{np\,+\,2p\,-n}}(Q_{R_{0}})}^{\frac{n+2}{np\,+\,p\,-n}}\right)^{\frac{n}{n\,+\,2}},
\end{split}
\end{equation}which holds for every $\varepsilon\in(0,1]$ and every sufficiently
small $h\in\mathbb{R}\setminus\left\{ 0\right\} $. Therefore, recalling
(\ref{eq:molli}) and letting $\varepsilon\rightarrow0$ in \eqref{eq:ee19},
by virtue of (\ref{eq:strong}) we obtain\begin{equation}\label{eq:ending}
\begin{split}
\int_{Q_{\varrho/2}(z_{0})}\left|\Delta_{s,h}H_{\frac{p}{2}}(Du)\right|^{2}dz\,&\leq\,c_{4}(n,p,R_{0})\left(\nu\,\Vert Df\Vert_{L^{\frac{np\,+\,4}{np\,+\,4\,-\,n}}(Q_{R_{0}})}+\,\Vert Df\Vert_{L^{\frac{np\,+\,4}{np\,+\,4\,-\,n}}(Q_{R_{0}})}^{\frac{np\,+\,4}{np\,+\,2\,-\,n}}\right)\\
&\,\,\,\,\,\,\,+\frac{c_{4}(n,p,R_{0})}{R^{2}}\left(\Vert Du\Vert_{L^{p}(Q_{R_{0}})}^{p}+\,\Vert Du\Vert_{L^{p}(Q_{R_{0}})}^{2}+\nu^{p}+\nu^{2}\right),
\end{split}
\end{equation}which holds for every $s\in\left\{ 1,\ldots,n\right\} $ and every
sufficiently small $h\in\mathbb{R}\setminus\left\{ 0\right\} $. This
proves the desired result for $p>2$. Moreover, letting $h\rightarrow0$
in the above inequality, we also obtain estimate \eqref{eq:faible}.\\
$\hspace*{1em}$Finally, when $p=2$, arguing in a similar fashion
we can reach the same conclusions.\end{proof}
\end{singlespace}
\begin{singlespace}

\section{The time derivative: proof of Theorem \ref{thm:timeregularity} \label{sec:timereg}}
\end{singlespace}

\begin{singlespace}
\noindent $\hspace*{1em}$This section is devoted to the study of
the existence and regularity of the time derivative of the weak solutions
to equation (\ref{eq:1}), under the assumptions of Theorem \ref{thm:main4}.
Indeed, we are now in position to give the\vspace{0.2cm}

\noindent \begin{proof}[\bfseries{Proof of Theorem~\ref{thm:timeregularity}}]
We shall keep both the notation and the parabolic cylinders used for
the proof of Theorem \ref{thm:main4}. Let us first show that
\begin{equation}
H_{p-1}(Du)\,\in\,L_{loc}^{p'}\left(0,T;W_{loc}^{1,p'}\left(\Omega,\mathbb{R}^{n}\right)\right).\label{eq:co2}
\end{equation}
Observe that for $p=2$ we have $H_{p-1}(Du)=H_{\frac{p}{2}}(Du)$,
so that assertion (\ref{eq:co2}) immediately follows from Theorem
\ref{thm:main4} in this case. Therefore, from now on we will assume
that $p>2$.\\
$\hspace*{1em}$Let us notice that for every $\xi\in\mathbb{R}^{n}$
we have 
\[
H_{p-1}(\xi)=\mathcal{F}(H_{\frac{p}{2}}(\xi)),
\]
where $\mathcal{F}:\mathbb{R}^{n}\rightarrow\mathbb{R}^{n}$ is the
function defined by
\[
\mathcal{F}(\eta):=\vert\eta\vert^{\frac{p-2}{p}}\eta,
\]
which is locally Lipschitz continuous for $p>2$. Thus, the function
$H_{p-1}(Du)$ is weakly differentiable with respect to the $x$-variable
by virtue of the chain rule in Sobolev spaces. From the definitions
of $\mathcal{F}$ and $H_{\frac{p}{2}}$, it follows that 
\begin{equation}
\left|D_{\eta}\,\mathcal{F}(H_{\frac{p}{2}}(Du))\right|\,\leq\,c_{1}\,\vert H_{\frac{p}{2}}(Du)\vert^{\frac{p-2}{p}}\leq\,c_{1}\,\vert Du\vert^{\frac{p-2}{2}},\label{eq:estF}
\end{equation}
for some positive constant $c_{1}\equiv c(n,p)$. Now, applying the
chain rule, the Cauchy-Schwarz inequality and estimate (\ref{eq:estF}),
we obtain\begin{equation}\label{eq:new001}
\begin{split}
\left|D_{x}\,H_{p-1}(Du)\right|^{p'}&\leq\,c_{2}(n,p)\,\vert D_{\eta}\,\mathcal{F}(H_{\frac{p}{2}}(Du))\vert^{p'}\,\vert D_{x}\,H_{\frac{p}{2}}(Du)\vert^{p'}\\
&\leq\,c_{3}\,\vert Du\vert^{\frac{(p-2)p'}{2}}\,\vert D_{x}\,H_{\frac{p}{2}}(Du)\vert^{p'},
\end{split}
\end{equation}where $c_{3}\equiv c_{3}(n,p)>0$. Using \eqref{eq:new001}, Hölder's
inequality with exponents $\left(\frac{2(p-1)}{p-2},\frac{2}{p'}\right)$
and estimate \eqref{eq:faible}, we get\begin{equation}\label{eq:new002}
\begin{split}
&\left(\int_{Q_{\varrho/2}(z_{0})}\vert DH_{p-1}(Du)\vert^{p'}dz\right)^{\frac{1}{p'}}\\
&\,\,\,\,\,\,\,\leq\,c_{4}(n,p)\,\Vert Du\Vert_{L^{p}(Q_{\varrho/2})}^{\frac{p-2}{2}}\,\left(\int_{Q_{\varrho/2}(z_{0})}\vert DH_{\frac{p}{2}}(Du)\vert^{2}\,dz\right)^{\frac{1}{2}}\\
&\,\,\,\,\,\,\,\leq\,c_{5}\,\,\Vert Du\Vert_{L^{p}(Q_{R_{0}})}^{\frac{p-2}{2}}\,\left[\nu\,\Vert Df\Vert_{L^{\vartheta}(Q_{R_{0}})}+\,\Vert Df\Vert_{L^{\vartheta}(Q_{R_{0}})}^{\frac{np\,+\,4}{np\,+\,2\,-\,n}}\right]^{\frac{1}{2}}\\
&\,\,\,\,\,\,\,\,\,\,\,\,\,\,+\,\,\frac{c_{5}}{R}\left[\Vert Du\Vert_{L^{p}(Q_{R_{0}})}^{2p-2}+\,\Vert Du\Vert_{L^{p}(Q_{R_{0}})}^{p}+(\nu^{p}+\nu^{2})\,\Vert Du\Vert_{L^{p}(Q_{R_{0}})}^{p-2}\right]^{\frac{1}{2}},
\end{split}
\end{equation}where $c_{5}$ is a positive constant depending on $n$, $p$, $\vartheta$
and $R_{0}$. Note that the right-hand side of \eqref{eq:new002}
is finite, and this implies (\ref{eq:co2}).\\
$\hspace*{1em}$Now, let $\kappa=\min\,\{\vartheta,p'\}$. Going back
to the weak formulation (\ref{eq:weaksol}), thanks to (\ref{eq:co2})
we can perform a partial integration in the second term on the left-hand
side with respect to the spatial variables. We thus obtain
\[
\int_{Q_{\varrho/2}(z_{0})}u\cdot\varphi_{t}\,\,dz\,=\,-\int_{Q_{\varrho/2}(z_{0})}\left(\sum_{\alpha=1}^{n}D_{\alpha}\left[\left(H_{p-1}(Du)\right)_{\alpha}\right]+f\right)\cdot\varphi\,\,dz,
\]
\\
for any $\varphi\in C_{0}^{\infty}(Q_{\varrho/2}(z_{0}))$, and the
desired conclusion immediately follows from \eqref{eq:faible} if
$p=2$, and from \eqref{eq:new002} if $p>2$, since $f\in L^{\vartheta}\left(0,T;W^{1,\vartheta}(\Omega)\right)$
with $\frac{np\,+\,4}{np\,+\,4\,-\,n}\leq\vartheta<\infty$ and
\[
\frac{np\,+\,4}{np\,+4-\,n}=\left(p+\frac{4}{n}\right)'<\,p'\,\,\,\,\,\,\,\,\,\,\,\mathrm{for\,\,every\,\,}p\geq2.
\]
$\hspace*{1em}$Furthermore, we can now observe that
\[
\partial_{t}u\,=\,\sum_{\alpha=1}^{n}D_{\alpha}\left[\left(H_{p-1}(Du)\right)_{\alpha}\right]+f\,\,\,\,\,\,\,\,\,\mathrm{in}\,\,\,\,Q_{\varrho/2}(z_{0}),
\]
from which we can infer\begin{equation}\label{eq:new003}
\begin{split}
\left(\int_{Q_{\varrho/2}(z_{0})}\left|\partial_{t}u\right|^{\kappa}\,dz\right)^{\frac{1}{\kappa}}&\leq\,\,n\,\Vert DH_{p-1}(Du)\Vert_{L^{\kappa}(Q_{\varrho/2}(z_{0}))}+\,\Vert f\Vert_{L^{\kappa}(Q_{\varrho/2}(z_{0}))}\\
&\leq\,\,c_{6}\,\Vert DH_{p-1}(Du)\Vert_{L^{p'}(Q_{\varrho/2}(z_{0}))}+\,c_{6}\,\Vert f\Vert_{L^{\vartheta}(Q_{\varrho/2}(z_{0}))},
\end{split}
\end{equation}where $c_{6}$ is a positive constant depending on $n$, $p$, $\vartheta$
and $R_{0}$. Thus, we can now deduce estimate \eqref{eq:estcor}
from \eqref{eq:new003} and \eqref{eq:faible} if $p=2$, and by combining
\eqref{eq:new002} and \eqref{eq:new003} when $p>2$.\end{proof}

\noindent \bigskip{}

\noindent $\hspace*{1em}$\textbf{Acknowledgements.} We gratefully
acknowledge Lorenzo Brasco for pointing out to us the reference \cite{Akh}.
Moreover, we would like to thank the reviewers for their valuable
comments, which helped to improve this work.
\end{singlespace}

\begin{singlespace}

\lyxaddress{\noindent \textbf{$\quad$}\\
$\hspace*{1em}$\textbf{Pasquale Ambrosio}\\
Dipartimento di Matematica e Applicazioni ``R. Caccioppoli''\\
Università degli Studi di Napoli ``Federico II''\\
Via Cintia, 80126 Napoli, Italy.\\
\textit{E-mail address}: pasquale.ambrosio2@unina.it}

\lyxaddress{\noindent $\hspace*{1em}$\textbf{Antonia Passarelli di Napoli}\\
Dipartimento di Matematica e Applicazioni ``R. Caccioppoli''\\
Università degli Studi di Napoli ``Federico II''\\
Via Cintia, 80126 Napoli, Italy.\\
\textit{E-mail address}: antpassa@unina.it}
\end{singlespace}

\end{document}